\documentclass[3p,preprint,12pt]{elsarticle}
\makeatletter\if@twocolumn\PassOptionsToPackage{switch}{lineno}\else\fi\makeatother

\usepackage[section]{placeins}
\usepackage{flafter}  
\usepackage{tabulary,xcolor}
\usepackage{amsfonts,amsmath,amssymb}
\usepackage[T1]{fontenc}
\usepackage{moreverb}
\usepackage{multicol}
\usepackage{amsmath}
\usepackage{relsize}
\usepackage{url,multirow,morefloats,floatflt,cancel,tfrupee}

\usepackage[capposition=top]{floatrow}
\usepackage{algorithm}
\usepackage[noend]{algpseudocode}
\usepackage{hyperref}
\usepackage{systeme}

\makeatletter
\let\save@ps@pprintTitle\ps@pprintTitle
\def\hlinewd#1{%
	\noalign{\ifnum0=`}\fi\hrule \@height #1%
	\futurelet\reserved@a\@xhline}

\AtBeginDocument{\ifNAT@numbers \biboptions{sort&compress}\fi}
\makeatother

\usepackage{etoolbox}
\makeatletter
\patchcmd{\hdots@for}{\hfill}{\hskip\z@\@plus 1filll}{}{}
\makeatother  

\usepackage{ifluatex}
\ifluatex
\usepackage{fontspec}
\defaultfontfeatures{Ligatures=TeX}
\usepackage[]{unicode-math}
\unimathsetup{math-style=TeX}
\else 
\usepackage[utf8]{inputenc}
\fi 
\ifluatex\else\usepackage{stmaryrd}\fi

\usepackage{url,multirow,morefloats,floatflt,cancel,tfrupee}
\makeatletter

\usepackage{indentfirst}
\AtBeginDocument{\@ifpackageloaded{textcomp}{}{\usepackage{textcomp}}}
\makeatother
\usepackage{colortbl}
\usepackage{xcolor}
\usepackage{pifont}
\usepackage[nointegrals]{wasysym}
\makeatletter

\def\mcWidth#1{\csname TY@F#1\endcsname+\tabcolsep}

\def\cAlignHack{\rightskip\@flushglue\leftskip\@flushglue\parindent\z@\parfillskip\z@skip}
\def\rAlignHack{\rightskip\z@skip\leftskip\@flushglue \parindent\z@\parfillskip\z@skip}

\if@twocolumn\usepackage{dblfloatfix}\fi 
\AtBeginDocument{
	\expandafter\ifx\csname eqalign\endcsname\relax
	\def\eqalign#1{\null\vcenter{\def\\{\cr}\openup\jot\m@th
			\ialign{\strut$\displaystyle{##}$\hfil&$\displaystyle{{}##}$\hfil
				\crcr#1\crcr}}\,}
	\fi
}

\AtBeginDocument{%
	\@ifpackageloaded{endfloat}%
	{\renewcommand\efloat@iwrite[1]{\immediate\expandafter\protected@write\csname efloat@post#1\endcsname{}}}{}%
}%

\let\lt=<
\let\gt=>
\def\processVert{\ifmmode|\else\textbar\fi}

\@ifundefined{subparagraph}{
	\def\subparagraph{\@startsection{paragraph}{5}{2\parindent}{0ex plus 0.1ex minus 0.1ex}%
		{0ex}{\normalfont\small\itshape}}%
}{}

\newcommand\role[1]{\unskip}
\newcommand\aucollab[1]{\unskip}

\@ifundefined{tsGraphicsScaleX}{\gdef\tsGraphicsScaleX{1}}{}
\@ifundefined{tsGraphicsScaleY}{\gdef\tsGraphicsScaleY{.9}}{}
\def\checkGraphicsWidth{\ifdim\Gin@nat@width>\linewidth
	\tsGraphicsScaleX\linewidth\else\Gin@nat@width\fi}

\def\checkGraphicsHeight{\ifdim\Gin@nat@height>.9\textheight
	\tsGraphicsScaleY\textheight\else\Gin@nat@height\fi}

\def\fixFloatSize#1{}
\let\ts@includegraphics\includegraphics

\def\inlinegraphic[#1]#2{{\edef\@tempa{#1}\edef\baseline@shift{\ifx\@tempa\@empty0\else#1\fi}\edef\tempZ{\the\numexpr(\numexpr(\baseline@shift*\f@size/100))}\protect\raisebox{\tempZ pt}{\ts@includegraphics{#2}}}}

\AtBeginDocument{\def\includegraphics{\@ifnextchar[{\ts@includegraphics}{\ts@includegraphics[width=\checkGraphicsWidth,height=\checkGraphicsHeight,keepaspectratio]}}}

\def\URL#1#2{\@ifundefined{href}{#2}{\href{#1}{#2}}}

\def\UrlOrds{\do\*\do\-\do\~\do\'\do\"\do\-}%
\g@addto@macro{\UrlBreaks}{\UrlOrds}
\makeatother

\usepackage{graphicx}
\usepackage{subcaption}
\graphicspath{ {Figures/} }
\begin{document}
	\begin{frontmatter}
		
		\title{Robust Data-Driven Discovery of Partial Differential Equations under Uncertainties}
		
		\author{Zhiming Zhang}
		\ead{zzhan506@asu.edu}  
		\author{Yongming Liu}
		\ead{corresponding author:yongming.liu@asu.edu}
		\address {School for Engineering of Matter, Transport and Energy, Arizona State University, Tempe, AZ, 85281, USA} 

\begin{abstract}
	Robust physics (e.g., governing equations and laws) discovery is of great interest for many engineering fields and explainable machine learning. A critical challenge compared with general training is that the term and format of governing equations is not known as a prior. In addition, significant measurement noise and complex algorithm hyperparameter tuning usually reduces the robustness of existing methods. A robust data-driven method is proposed in this study for identifying the governing Partial Differential Equations (PDEs) of a given system from noisy data. The proposed method is based on the concept of Progressive Sparse Identification of PDEs (PSI-PDE or $\psi$-PDE). Special focus is on the handling of data with huge uncertainties (e.g., 50$\%$ noise level). Neural Network modeling and fast Fourier transform (FFT) are implemented to reduce the influence of noise in sparse regression. Following this, candidate terms from the prescribed library are progressively selected and added to the learned PDEs, which automatically promotes parsimony with respect to the number of terms in PDEs as well as their complexity. Next, the significance of each learned terms is further evaluated and the coefficients of PDE terms are optimized by minimizing the L2 residuals. Results of numerical case studies indicate that the governing PDEs of many canonical dynamical systems can be correctly identified using the proposed $\psi$-PDE method with highly noisy data. One great benifit of proposed algorithm is that it avoids complex algorithm modification and hyperparameter tuning in most existing methods. Limitations of the proposed method and major findings are presented. 
\end{abstract}

\end{frontmatter}
\section{Introduction}\label{Intro}
Despite that many dynamical systems can be well characterized by PDEs derived mathematically/physically from basic principles such as conservation laws, lots of other systems have unclear or elusive underlying mechanisms (e.g., ones in neuroscience, finance, and ecology). Thus, the governing equations are usually empirically formulated \cite{rudy2017data}. Data-driven physics discovery of dynamical systems gradually became possible in recent years due to the rapid development and extensive application of sensing technologies and computational power \cite{long2019pde}. Over the past years, extensive efforts have been devoted into discovering representative PDEs for complex dynamical systems of which limited prior knowledge are available \cite{schmidt2009distilling,raissi2018hidden,rudy2017data,long2019pde}. 

Among all the methods investigated for PDE identification \cite{schmidt2009distilling,raissi2018hidden,rudy2017data,long2019pde,maslyaev1903data,atkinson2019data,hasan2020learning,xu2020dlga}, sparse regression gains the most attention in recent studies due to its inherent parsimony-promoting advantage. Considering a nonlinear PDE of the general form $u_t = N(u,u_x,u_{xx},...,x)$, in which the subscripts denote partial differentiation with respect to temporal or spatial coordinate(s), $N(\cdot)$ is an unknown expression on the right hand side of the PDE. It is usually a nonlinear function of the spatial coordinate $x$, the measured quantity $u(x,t)$, and its spatial derivatives $u_x$,$u_{xx}$, etc. Given time series measurements of $u$ at certain spatial locations, the above equation can be approximated as $\mathbf{U}_t=\boldsymbol{\Theta}(\mathbf{U})\xi$, in which $\mathbf{U}_t$ is the discretized form of $u_t$, $\boldsymbol{\Theta}(\mathbf{U})$ is a library matrix with each column corresponding to a candidate term in $N(\cdot)$. A key assumption in sparse identification is that $N(\cdot)$ consists of only a few term for a real physical system, which requires the solution of regression (i.e., $\xi$) to be a sparse vector with only a limited number of nonzero elements. This assumption promotes a parsimonious form of the learned PDE instead of overfitting the measured data with a complex model containing redundant nonlinear higher-order terms.

As pioneering researchers in sparse PDE learning, Rudy et al. \cite{rudy2017data,rudy2019data} modified the ridge regression method by imposing hard thresholding which recursively eliminates certain terms with coefficient values below a learned threshold. As pointed out in Limitations of \cite{rudy2017data,rudy2019data} (Section 4 in Supplementary Materials) and following studies \cite{raissi2018hidden,chen2020deep,both2020sparsely}, the identification quality is very sensitive to data quantity and quality. For example, the terms of the reaction diffusion equation cannot be correctly identified using the data with only 0.5\% random noise. Furthermore, as indicated in \cite{fuentes2020equation}, the identification results using this method are susceptible to the selection of hyperparameters of the algorithm, including the regularizer $\lambda$ and the initial tolerance which is also the tolerance increment $d_{tol}$. The hyperparameter tuning is especially critical for cases with noisy measurements. This limitation most probably comes from the hard thresholding in the modified algorithm (STRidge). A hard thresholding tends to suppress small coefficients that may not correspond to the most trivial terms of the intermediately learned PDEs. 

To overcome the challenge of numerical differentiation with scarce and noisy data in sparse regression methods, deep learning techniques were incorporated by generating a large quantity of meta-data and adopting the automatic differentiation function in deep learning frameworks (Tensorflow, PyTorch, etc.) \cite{xu2019dl,berg2019data}. The intermediately learned PDE can be treated as a physics loss term in physics-informed deep learning \cite{raissi2019physics,wang2020towards,haghighat2020sciann}, and constrained neural networks were developed to improve the performance of PDE identification recursively \cite{hasan2020learning,chen2020deep,both2020sparsely}. Long et al. \cite{long2018pde,long2019pde} used a convolutional architecture and symbolic regression to replace the numerical differentiation and sparse regression procedures, respectively. A comprehensive review of the state of the art of PDE learning can be found in \cite{chen2020deep}. Despite improved performance of PDE identification using these methods, the identification results (both PDE forms and coefficients) are lacking robustness in most studies mentioned above. For example, approaches using constrained neural networks introduced more hyperparameters into the algorithms in addition to those in the used STRidge algorithm, which further increases the challenge of identifying the correct PDE forms since PDE learning problems are sensitive to the hyperparameter tuning. This issue is amplified under noisy data, especially under high noise levels. Complete different identification results may be obtained under different noise levels using same hyperparameter settings. Thus, a sound robust PDE learning needs to produce stable identification results with respect to different noise levels. 

Considering the gaps of existing studies in discovering PDEs from complex dynamical systems, a robust method for correctly identifying PDEs is needed to discover the underlying physics of the measured systems that lack prior knowledge of the governing principles. Thus, this study attempts to develop a robust method of PDE identification within the framework of sparse regression. The key idea is to address both sparsity and accuracy of the learned PDE. Special focus is on the automatic and progressive selection of learned PDE forms without complex algorithms with hard-to-tune hyperparameters \cite{schmidt2009distilling}.  The proposed scheme automatically promotes sparsity in addition to simplicity of the learned model. Finally, the representativeness of each model will be further evaluated by solving its corresponding PDE with given/extracted initial and/or boundary conditions. The coefficients of each term are optimized by minimizing the error of model prediction with the measured data taken as the ground truth. In this way, the PDE that is most likely to represent the intrinsic mechanisms underlying the observed system will be determined. Since the proposed methodology progressively yields a sparse identification of the governing PDE(s) of a given system, it is named the progressive and sparse identification method of PDEs (PSI-PDE or $\psi$-PDE method).

The remaining part of this paper is structured as follows. Section \ref{Sec:method} establishes the framework of the $\psi$-PDE method; section \ref{Sec:results} presents and discusses the results of discovering govern equations for a variety of dynamical systems using the $\psi$-PDE method; section \ref{Section:Conclusion} concludes this study with remarks and recommendations for future work.

\section{Methodology: a robust PDE learning method}\label{Sec:method}

\begin{figure}[!h]
	\centering
	\includegraphics[scale=0.5]{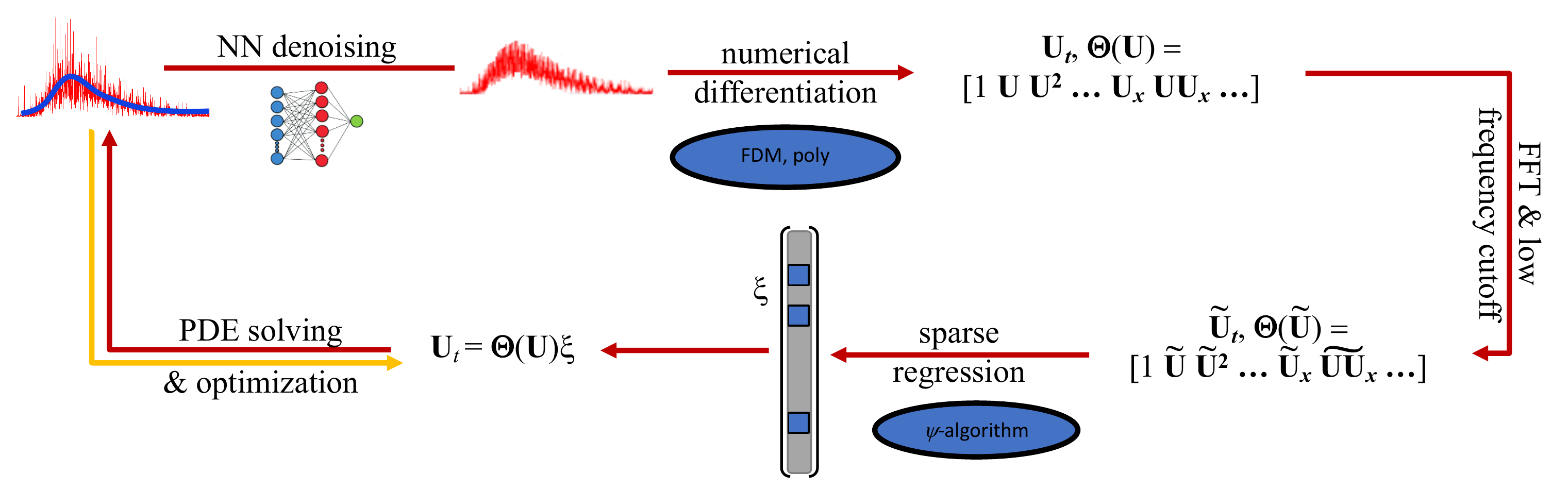}
	\caption{Framework of the $\psi$-PDE method for discovering PDE(s) from measured data.}
	\label{Figure:Framework}
\end{figure}

Figure \ref{Figure:Framework} illustrates the framework of the proposed $\psi$-PDE method for discovering the governing equations using measured data from a dynamical system. This framework starts from the red noisy curve on the upper left corner, which denotes the measured signals containing all the information one can directly obtain from the instrumented system. For preprocessing the measured data, a neural network (NN) model is built following the practices in \cite{xu2019dl,berg2019data}, setting the independent variables (i.e., $t$, $x$, etc.) as inputs and the measured quantity (e.g., $u$) as the output. The measured data are split into training and validation sets and the early stopping strategy is devised in the model training to prevent the NN model from overfitting the noise components contained in the measured data. A smoothed series of signals is expected from this preprocessing, which will be subsequently used to calculate the numerical derivatives (i.e., $\mathbf{U}_t$, $\mathbf{U}_x$, $\mathbf{U}_{xx}$, etc.) and then construct the library matrix $\boldsymbol{\Theta}(\mathbf{U})$ for sparse regression. Numerical methods such as finite difference method (FMD) and polynomial interpolation are used to calculate the temporal and spatial derivatives. With $\mathbf{U}_t$ and $\boldsymbol{\Theta}(\mathbf{U})$ established, to further reduce the influence of noise in sparse regression, fast Fourier transform (FFT) is applied to transform $\mathbf{U}_t$ and $\boldsymbol{\Theta}(\mathbf{U})$ to their frequency domain counterparts $\widetilde{\mathbf{U}}_t$ and $\boldsymbol{\Theta}(\widetilde{\mathbf{U}})$, respectively. Following FFT, a frequency cutoff is implemented to preserve only the low frequency components that are expected to be less susceptible to noise. Moreover, this step converts the regression problem from the temporal-spatial domain to the frequency domain, which does not change the form of learned PDEs \cite{cao2020machine}. 

With $\widetilde{\mathbf{U}}_t$ and $\boldsymbol{\Theta}(\widetilde{\mathbf{U}})$ from FFT with frequency cutoff, sparse regression is conducted in a progressive manner. Algorithms that automatically yield a sparse solution by imposing $\ell_1$ norm regularization (such as LASSO) or hard thresholding (such as STRidge) are not adopted to avoid the lemma of hyperparameter tuning in PDE learning. Instead, a least squares regression is implemented via $x = A\backslash b$ in MATLAB to recursively examine the importance of each term in the prescribed library by evaluating the resulting regression error and model complexity. In this way, the most important term(s) are step-by-step identified and added to the PDE model until the effects of adding more terms diminish. The details of this procedure are elaborated in the $\psi$-algorithm in Algorithm \ref{alg1}. This algorithm probably yields more than one candidate PDE models that are hard to compare from the perspective of regression accuracy and model complexity. Finally, all candidate PDEs are solved numerically given sufficient initial/boundary conditions, and the solutions (the blue smooth curve at the upper left corner of Figure \ref{Figure:Framework}) are compared with the measured data. In this step, the importance of each term can be further evaluated by eliminating certain terms and optimizing their coefficients. The final PDE for the given system is determined as the one capturing the most intrinsic mechanism represented by the essential terms. The rest of this section demonstrates each step of the $\psi$-PDE method using an example of discovering the Burgers equation from simulated data. 

\newcommand{\norm}[1]{\left\lVert#1\right\rVert}
\begin{algorithm}
	\caption{$\psi$-algorithm: $\xi = \psi(\boldsymbol{\Theta},\widetilde{\mathbf{U}}_t,\gamma_{Reg},\gamma_{BIC})$}\label{alg1}
	\begin{algorithmic}[1]
		
		\State \textbf{Input}: library matrix $\boldsymbol{\Theta}$ (with the size $M\times N$), discretized temporal derivative $\widetilde{\mathbf{U}}_t$ (with the size $M\times 1$), tolerance ratio of root-mean-square regression error $\gamma_{Reg}$, and tolerance ratio of BIC $\gamma_{BIC}$.
		\State Normalize $\boldsymbol{\Theta}$ and $\widetilde{\mathbf{U}}_t$:
		$\boldsymbol{\Theta} = \boldsymbol{\Theta}/\norm{\boldsymbol{\Theta}}_2$,
		$\widetilde{\mathbf{U}}_t = \widetilde{\mathbf{U}}_t/\norm{\widetilde{\mathbf{U}}_t}_2$
		
		\textcolor{gray}{\texttt{\#} Steps 3 to 7 identify the most contributive terms $I_0$.}
		\For{$i = 1,2,...,n_{val}$}	\qquad \textcolor{gray}{\texttt{\#} $n_{val} = 1\mathrm{E}4$: number of validations.}
		
			{
				\State Randomly split data with $n_{trn}/n_{val} = 80/20$:
				$\boldsymbol{\Theta}^{trn}$,
				$\boldsymbol{\Theta}^{val}$,
				$\widetilde{\mathbf{U}}_t^{trn}$, and
				$\widetilde{\mathbf{U}}_t^{val}$
				\For{j =  1:N} \qquad \textcolor{gray}{\texttt{\#} evaluate the importance of $j^\mathrm{th}$ term in $\boldsymbol{\Theta}$.}
				
				{
					\quad$indSel = setdiff(1:N,j)$	\qquad \textcolor{gray}{\texttt{\#} delete $j$ from the selected index list.}
				
					\quad$\boldsymbol{\Theta}_1 = \boldsymbol{\Theta}^{trn}(:,indSel)$
				
					\quad$\boldsymbol{\Theta}_2 = \boldsymbol{\Theta}^\mathrm{val}(:,indSel)$
				
					\quad$\xi_1 = \boldsymbol{\Theta}_1  \setminus \widetilde{\mathbf{U}}_t^{trn}$ 	\qquad \textcolor{gray}{\texttt{\#} $\setminus$: matrix left division.}
				
					\quad$\widetilde{\mathbf{U}}_{t2} = \boldsymbol{\Theta}_2 \times \xi_1$
				 
				 	\quad$\varepsilon_{Reg}(i,j) = rms(\widetilde{\mathbf{U}}_{t2}-\widetilde{\mathbf{U}}_t^{val})$
				 	\qquad \textcolor{gray}{\texttt{\#} root mean square of regression error.}
				 	
				 	\quad$\varepsilon_{BIC}(i,j) = n_{trn} * \mathrm{log}(mse(\widetilde{\mathbf{U}}_{t2},\widetilde{\mathbf{U}}_t^{val}))+\frac{\mathrm{sum}(indSel^2)+(N-1)^2}{N-1}\mathrm{log}(n_{trn})$
				 	
				 	\quad\textcolor{gray}{\texttt{\#} BIC, the term $\frac{\mathrm{sum}(indSel^2)}{N-1}$ is added to penalize the selected complex terms.}
				}
				\EndFor	
			}	
		\EndFor	
		\State Plot histograms  and mean values of $\varepsilon_{Reg}$ and $\varepsilon_{BIC}$ for each term in $\boldsymbol{\Theta}$.
		\State Identify the term(s) with index/indices ($I_0$) that correspond to the largest regression error ($\varepsilon_{Reg}$) and/or BIC ($\varepsilon_{BIC}$), considering both the distributions and mean values. 
	
		\textcolor{gray}{\texttt{\#} The following \textbf{while} loop progressively adds important terms to $I_0$.}
		\While{$std(mean(\varepsilon_{Reg}))>\gamma_{Reg}*\varepsilon_{Reg}^{ref}$ \textbf{or} $std(mean(\varepsilon_{BIC}))>\gamma_{BIC}*\varepsilon_{BIC}^{ref}$ }

			\For{$i = 1,2,...,n_{val}$}	\qquad \textcolor{gray}{\texttt{\#} $n_{val} = 1\mathrm{E}4$: number of validations.}
			
			{
				\State Randomly split data with $n_{trn}/n_{val} = 80/20$:
				$\boldsymbol{\Theta}^{trn}$,
				$\boldsymbol{\Theta}^{val}$,
				$\widetilde{\mathbf{U}}_t^{trn}$, and
				$\widetilde{\mathbf{U}}_t^{val}$
				\For{j =  1:N}
				
				{
					\qquad$indSel = \mathrm{union}(I_0,j)$	\qquad \textcolor{gray}{\texttt{\#} add $j$ to the selected index list.}
					
					\qquad$\boldsymbol{\Theta}_1 = \boldsymbol{\Theta}^{trn}(:,indSel)$
					
					\qquad$\boldsymbol{\Theta}_2 = \boldsymbol{\Theta}^\mathrm{val}(:,indSel)$
					
					\qquad$\xi_1 = \boldsymbol{\Theta}_1  \setminus \widetilde{\mathbf{U}}_t^{trn}$ 
					
					\qquad$\widetilde{\mathbf{U}}_{t2} = \boldsymbol{\Theta}_2 \times \xi_1$
					
					\qquad$\varepsilon_{Reg}(i,j) = rms(\widetilde{\mathbf{U}}_{t2}-\widetilde{\mathbf{U}}_t^{val})$
					
					\qquad$\varepsilon_{BIC}(i,j) = n_{trn} * \mathrm{log}(mse(\widetilde{\mathbf{U}}_{t2},\widetilde{\mathbf{U}}_t^{val}))+\frac{\mathrm{sum}(indSel^2)+(n_{I_0}+1)^2}{n_{I_0}+1}\mathrm{log}(n_{trn})$
					
				}
				\EndFor	
				$\varepsilon_{Reg}^{ref}(i) = rms(\widetilde{\mathbf{U}}_t^{val})$
				
				\quad $\varepsilon_{BIC}^{ref}(i) = n_{trn} * \mathrm{log}(mse(0,\widetilde{\mathbf{U}}_t^{val}))+\frac{\mathrm{sum}((1:N)^2)+N^2}{N}\mathrm{log}(n_{trn})$
			}	
			\EndFor	
			\State $\varepsilon_{Reg}^{ref} = mean(\varepsilon_{Reg}^{ref})$
			\State $\varepsilon_{Reg}^{BIC} = mean(\varepsilon_{Reg}^{BIC})$
			\State Plot histograms  and mean values of $\varepsilon_{Reg}$ and $\varepsilon_{BIC}$ for each selection of terms in $\boldsymbol{\Theta}$.
			\State Identify the term(s) with index/indices ($I_a$) that correspond to the smallest regression error ($\varepsilon_{Reg}$) and/or BIC ($\varepsilon_{BIC}$), considering both the distributions and mean values.
			\State $I_0 = \mathrm{union}(I_0,I_a)$
		\EndWhile


	\end{algorithmic}
\end{algorithm}

\begin{figure}[!h]
	\centering
	\includegraphics[scale=0.8]{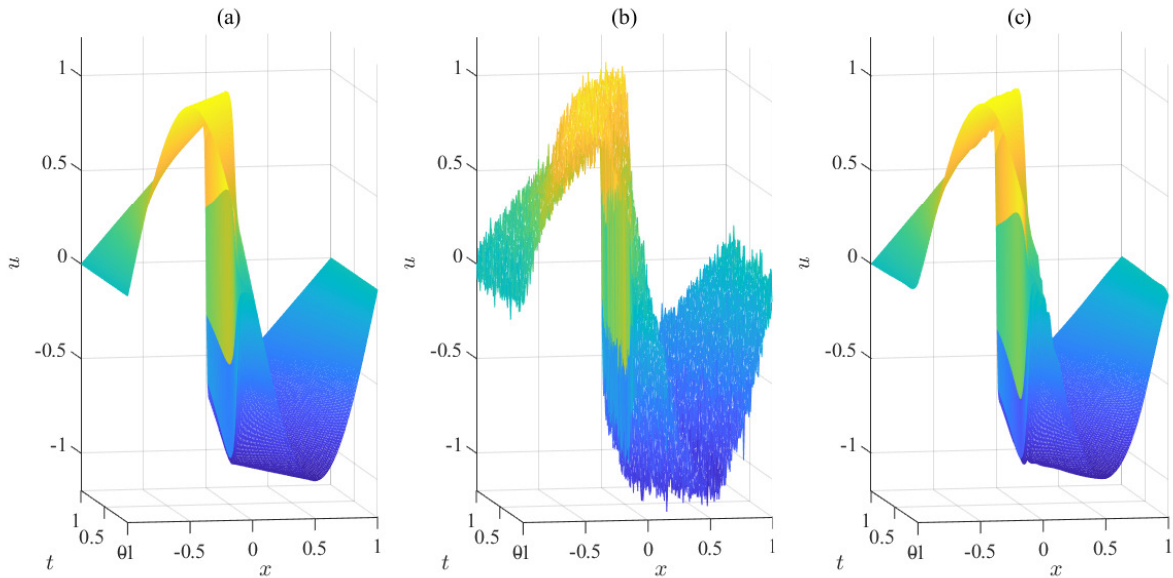}
	\caption{Solution of Burgers equation $u_t = -uu_x+\frac{0.01}{\pi} u_{xx}$ with (a) 0\% noise, (b) 10\% noise, and (c) 10\% noise afer NN denoising.}
	\label{Figure:solBurgers0}
\end{figure}

Burgers equation is used to describe the dynamics of a dissipative system. A 1D viscous Burgers equation is used to explain the steps of the $\psi$-PDE method in detail. It has the expression of $u_t = -uu_x+\nu u_{xx}$ with the initial and boundary conditions $u(0,x) = -\mathrm{sin}(\pi x)$ and $u(t,-1) = u(t,1) =0$, in which $\nu = \frac{0.01}{\pi}$ denotes the diffusion coefficient. Figure \ref{Figure:solBurgers0} (a) shows its solution within the range $t\in [0,1]$ and $x\in [-1,1]$. The library for sparse regression is built with polynomials of $u$ to the power of 3, spatial derivatives to the $3^\mathrm{rd}$ order, and their products. As a result, it contains 16 terms in total, i.e., $\boldsymbol{\Theta}(\mathbf{U}) = \{1, \mathbf{U}, \mathbf{U}^2, \mathbf{U}^3, \mathbf{U}_x, \mathbf{UU}_x,...,  \mathbf{U}^3\mathbf{U}_{xxx}\}$. 

To demonstrate the effect of NN denoising step in the $\psi$-PDE method, 10\% white Gaussian noise is added to the numerical solution of the Burgers equation, which significantly varies the values of the solution (as shown in Figure \ref{Figure:solBurgers0} (b)) and poses challenge to calculating the numerical derivatives. In this study, noise level is quantified by the percentage of the standard deviation of the measured variable. For example, if 10\% noise is added to $u$, then the outcome is $u_n = u+10\%*std(u)*randn(size(u))$ where $randn(\cdot)$ generates white Gaussian noise of the specified dimension. Without much prior knowledge about the noise characteristics, NN modeling is applied to denoise the noisy measurements, and the processed data is visualized in Figure \ref{Figure:solBurgers0} (c). Comparing the three plots in Figure \ref{Figure:solBurgers0}, one can observe that NN denoising largely reduces the noise level in the collected data (from 10\% to 2\%) and makes the solution curve much smoother than the noisy one, which has the potential of improving the accuracy of subsequent numerical differentiation.

\begin{figure}[!h]
	\centering
	\includegraphics[scale=0.8]{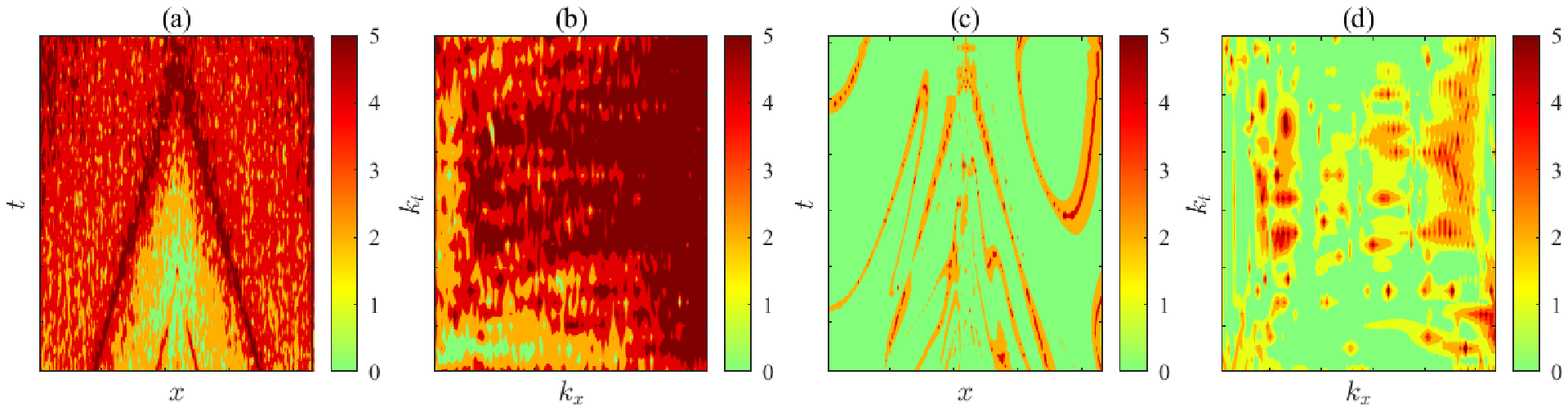}
	\caption{Color maps of relative errors in $u_{xx}$ of Burger equation caused by adding 50\% noise. 
		(a) $\mathrm{log}(\left|\frac{u^n_{xx}-u^0_{xx}}{u^0_{xx}}\right|)$;
		(b) $\mathrm{log}(\left|\frac{\widetilde{u}^n_{xx}-\widetilde{u}^0_{xx}}{\widetilde{u}^0_{xx}}\right|)$;
		(c) $\mathrm{log}(\left|\frac{u^{nn}_{xx}-u^0_{xx}}{u^0_{xx}}\right|)$;
		(d)$\mathrm{log}(\left|\frac{\widetilde{u}^{nn}_{xx}-\widetilde{u}^0_{xx}}{\widetilde{u}^0_{xx}}\right|)$. $u^0_{xx}$, $u^n_{xx}$, and $u^{nn}_{xx}$ are the $2^\mathrm{nd}$ order derivatives calculated using the clean data, noise data, and NN-denoised data, respectively; the tilde $\widetilde{(\cdot)}$ denote the results of 2D FFT; $k_x$ and $k_t$ represent the corresponding coordinates in the frequency domain.}
	\label{Figure:effectFFT}
\end{figure}

With numerical derivatives calculated and the library matrix established, FFT is conducted to further reduce the influence of noise in the subsequent sparse regression. Figures \ref{Figure:effectFFT} (a) to (d) compare the relative errors in the spatial-temporal domain and frequency domain after taking 2D FFT. Figures \ref{Figure:effectFFT} (a) and (b) show the difference between the polluted data with 50\% noise and the simulated clean data. It can be observed that after taking 2D FFT, the low-frequency components are less affected by the added noise. In addition, Figures \ref{Figure:effectFFT} (c) and (d) show that NN denoising can largely reduces not only the relative error in the spatial-temporal domain, as can be predicted from Figure \ref{Figure:solBurgers0}, but also the relative error in the frequency domain. Therefore, it can be expected that the performance of PDE learning can be considerably improved by implementing these preprocessing procedures.

\begin{figure}[!h]
	\centering
	\tiny{(a) \hspace{130pt} (b)  \hspace{130pt} (c)}\\
	\includegraphics[scale=0.4]{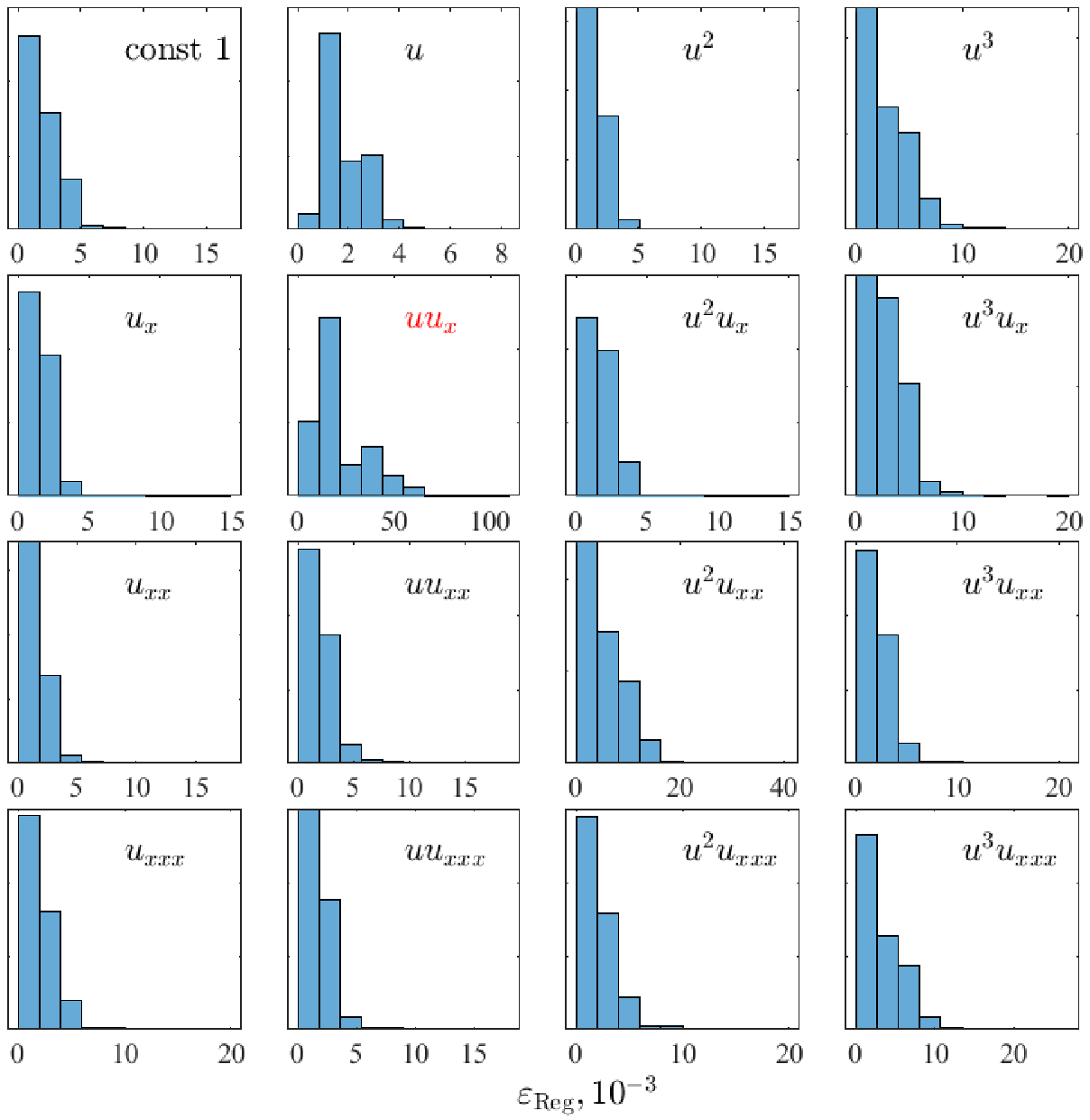}
	\includegraphics[scale=0.4]{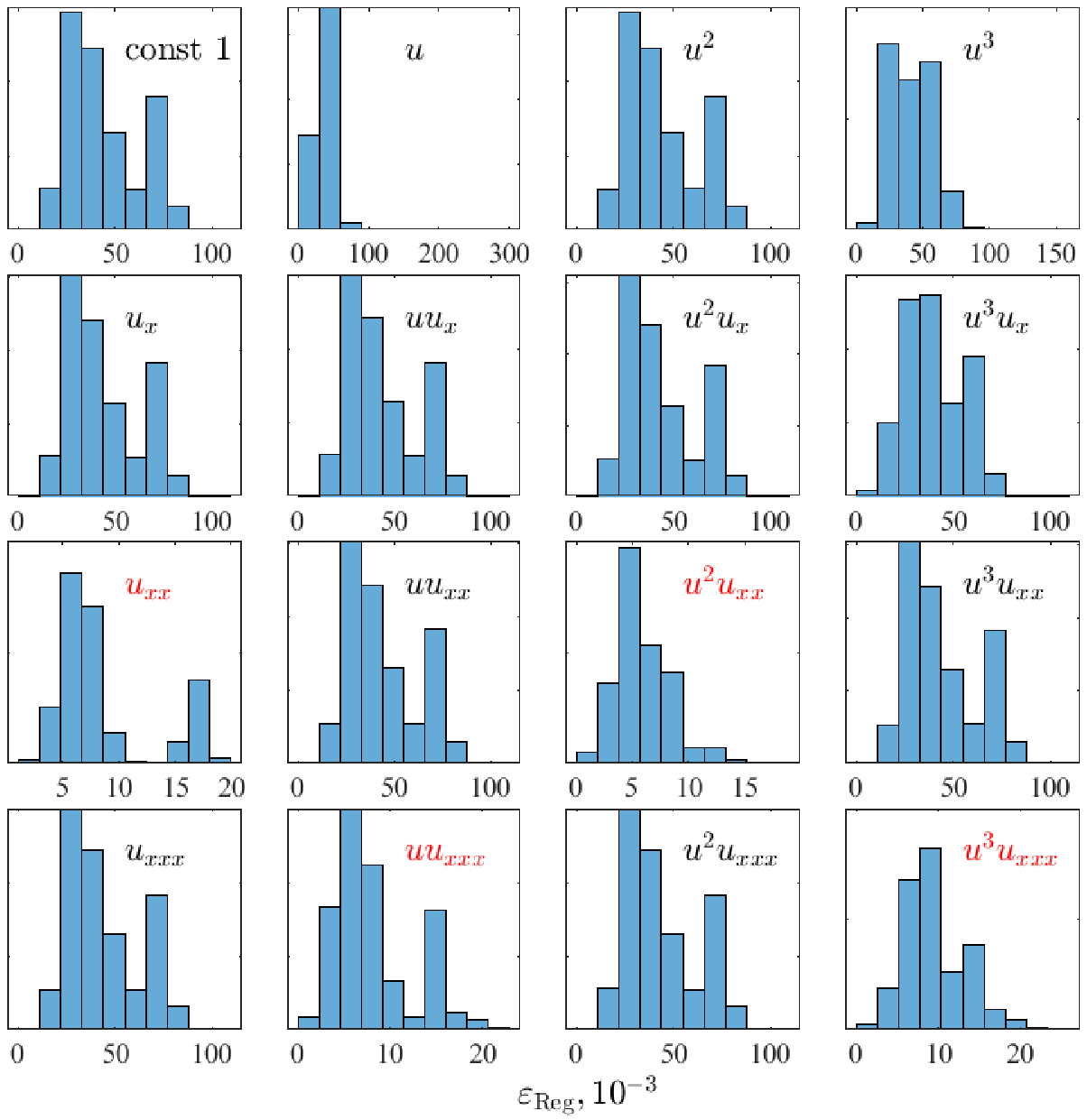}
	\includegraphics[scale=0.4]{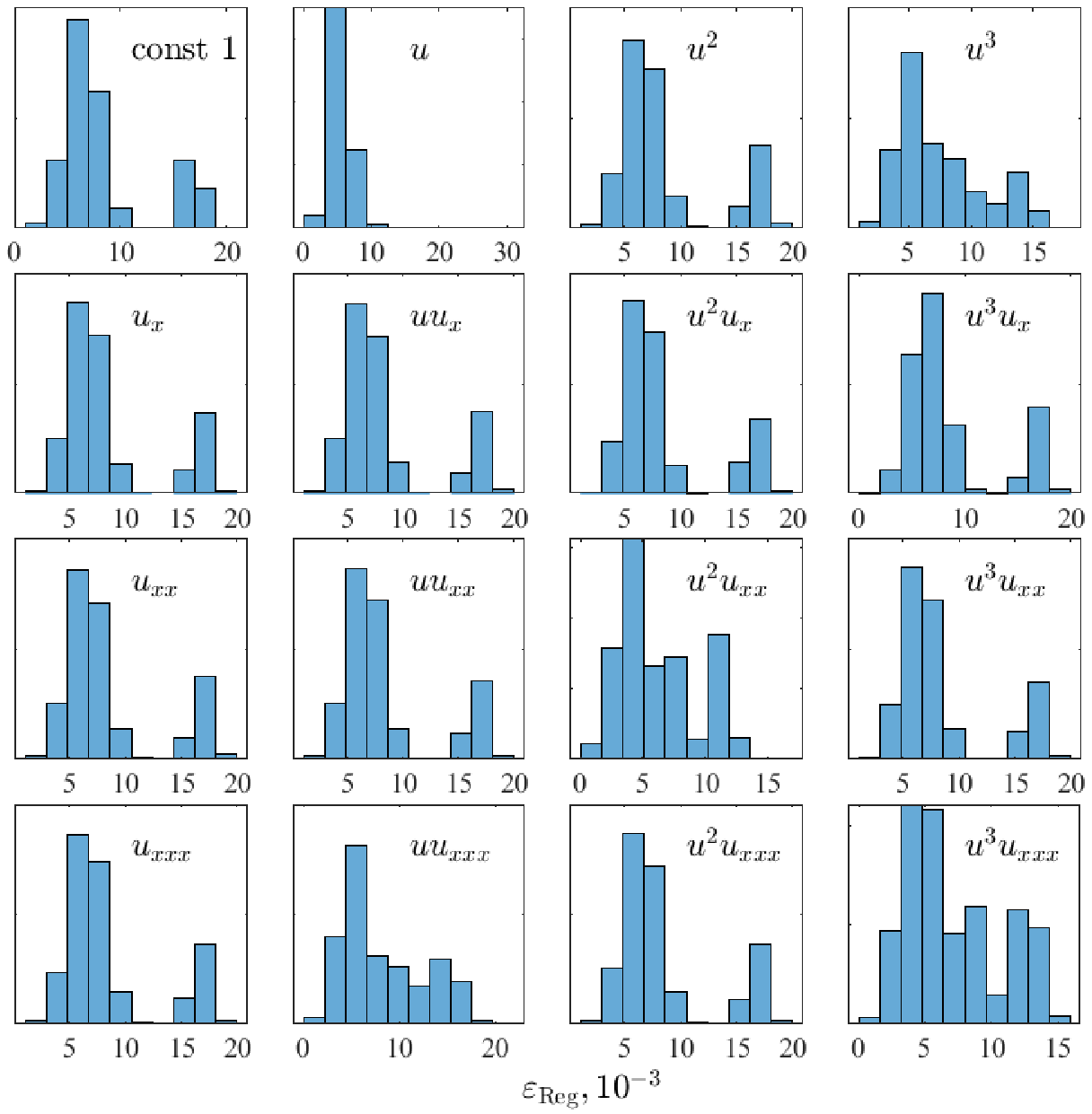}
	\caption{Distribution of regression error when dropping (a) or adding (b and c) a certain term in the library (Burgers equation: $u_t = -uu_x+\frac{0.01}{\pi}u_{xx}$).}
	\label{Figure:errReg}
\end{figure}

\begin{figure}[!h]
	\centering
	\tiny{(a) \hspace{130pt} (b)  \hspace{130pt} (c)}\\
	\includegraphics[scale=0.4]{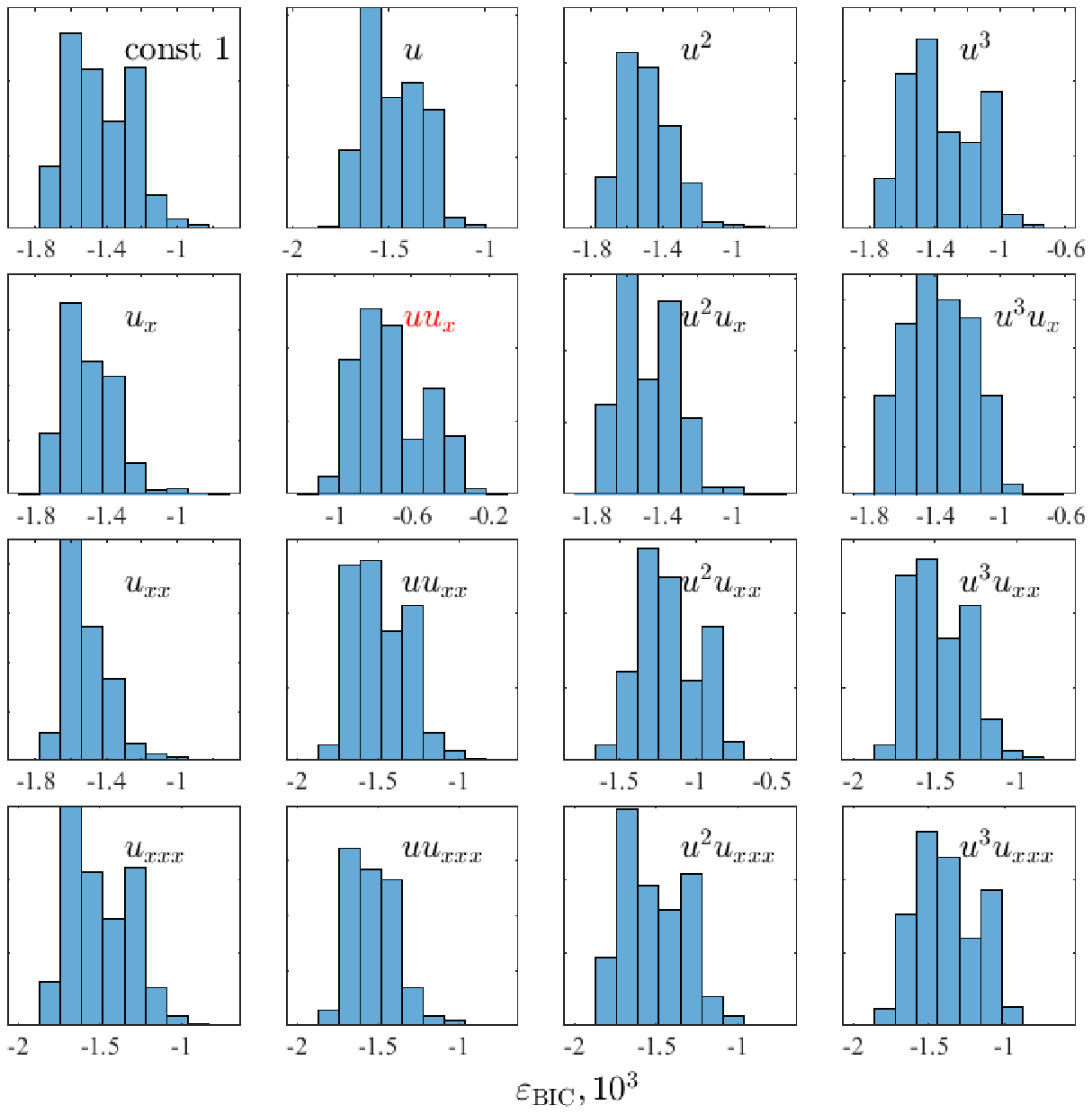}
	\includegraphics[scale=0.4]{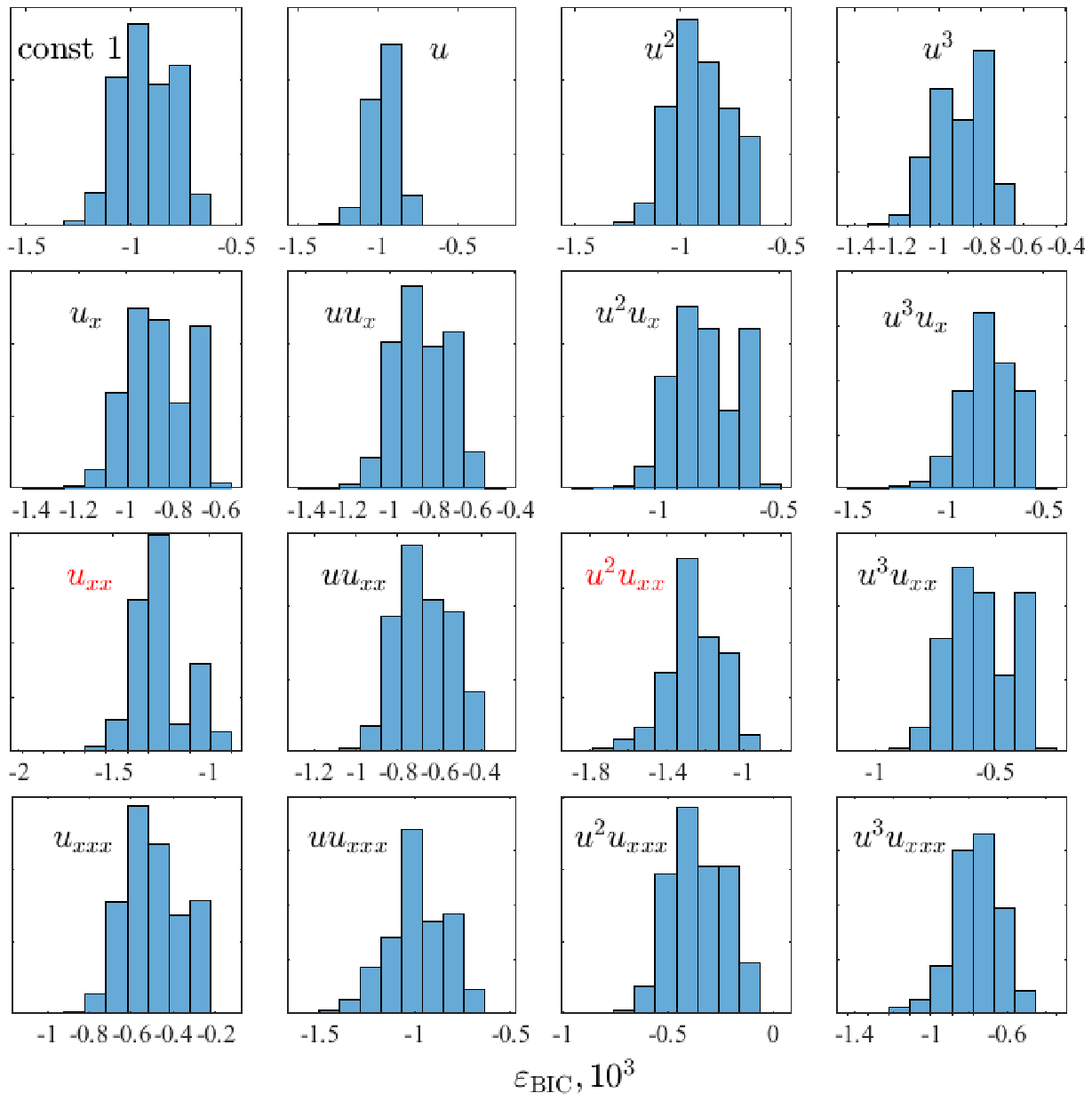}
	\includegraphics[scale=0.4]{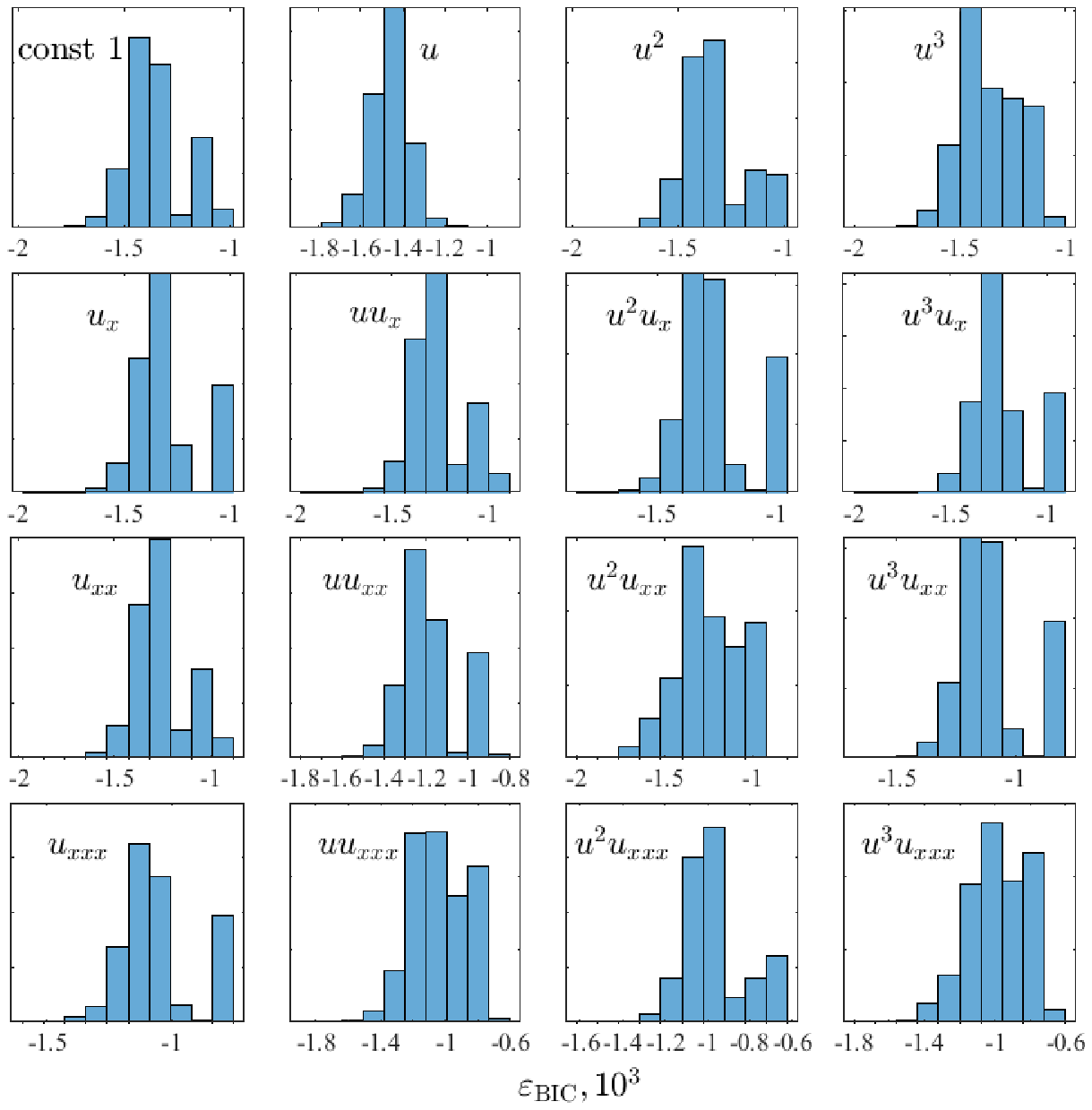}
	\caption{Distribution of BIC when dropping (a) or adding (b and c) a certain term in the library (Burgers equation: $u_t = -uu_x+\frac{0.01}{\pi}u_{xx}$).}
	\label{Figure:errBIC}
\end{figure}

\begin{figure}[!h]
	\centering
	\tiny{(a) \hspace{130pt} (b)  \hspace{130pt} (c)}\\
	\includegraphics[scale=0.4]{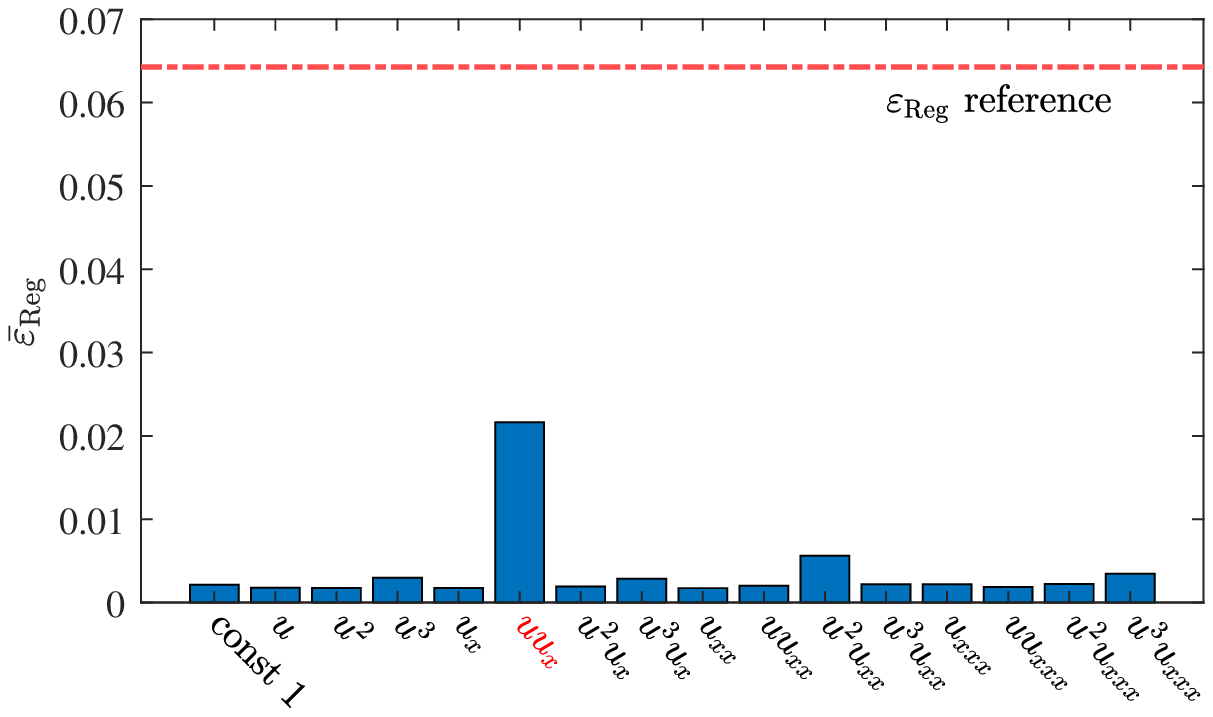}
	\includegraphics[scale=0.4]{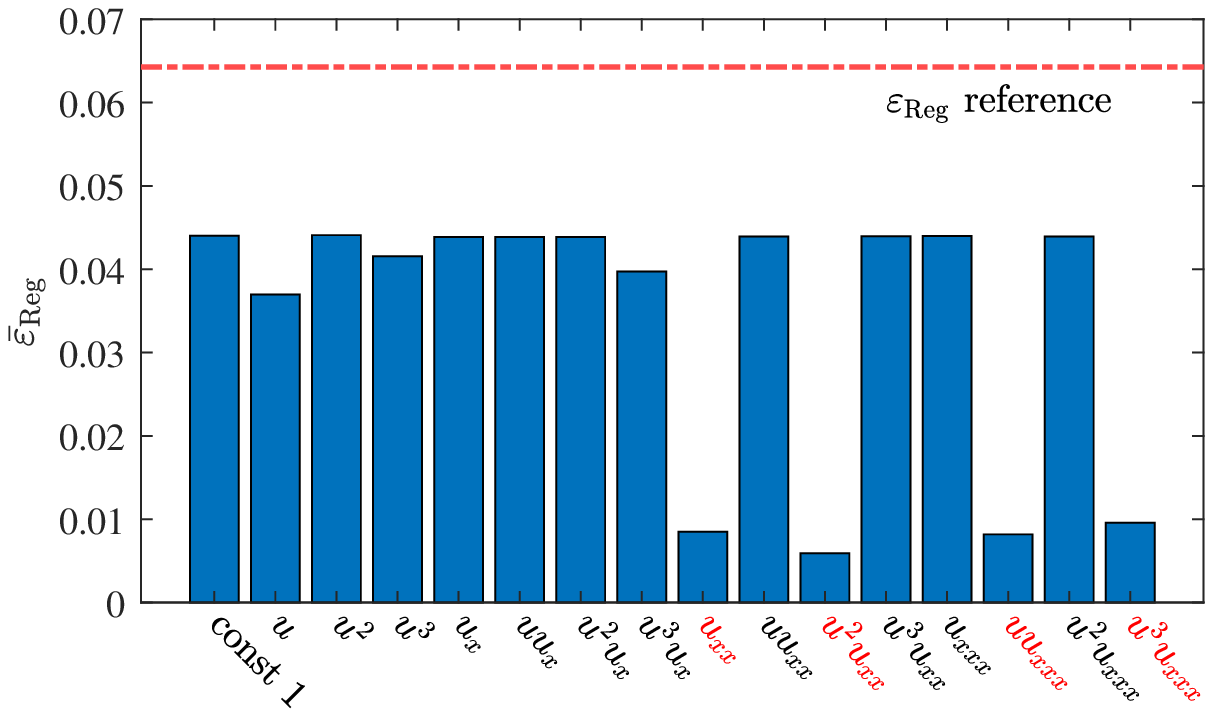}
	\includegraphics[scale=0.4]{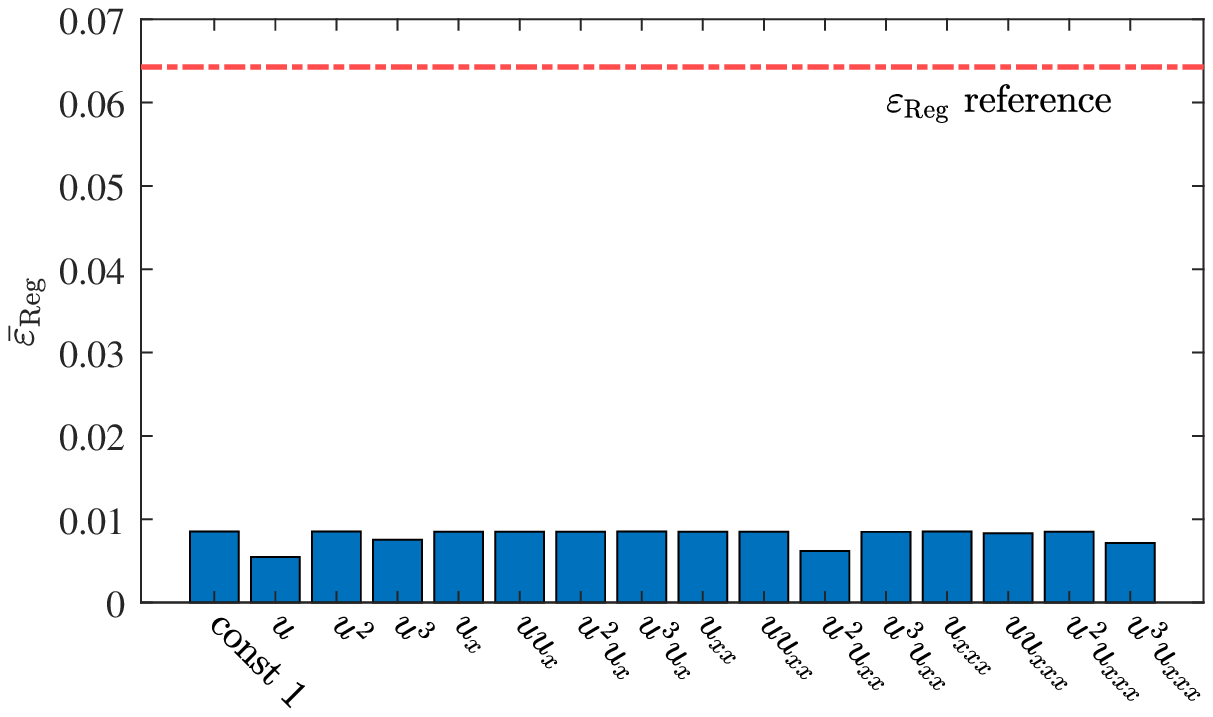}
	\caption{Mean regression error when dropping (a) or adding (b and c) a certain term in the library (Burgers equation: $u_t = -uu_x+\frac{0.01}{\pi}u_{xx}$).}
	\label{Figure:aveReg}
\end{figure}

\begin{figure}[!h]
	\centering
	\tiny{(a) \hspace{130pt} (b)  \hspace{130pt} (c)}\\
	\includegraphics[scale=0.4]{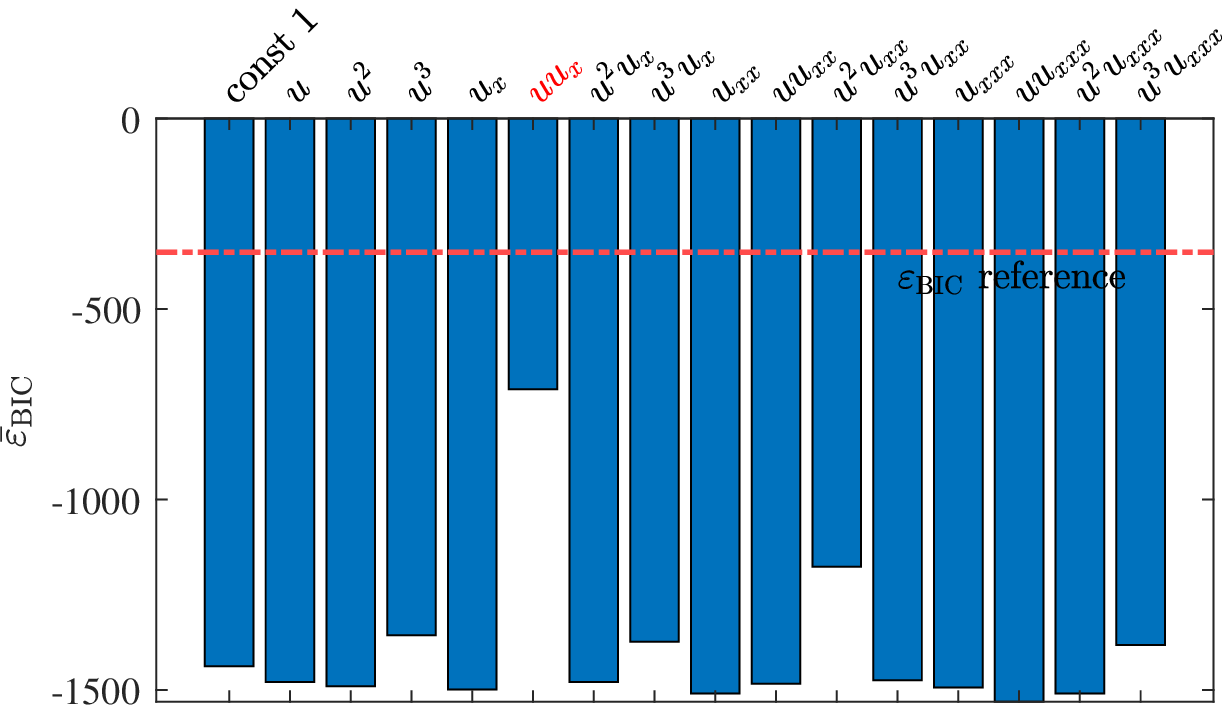}
	\includegraphics[scale=0.4]{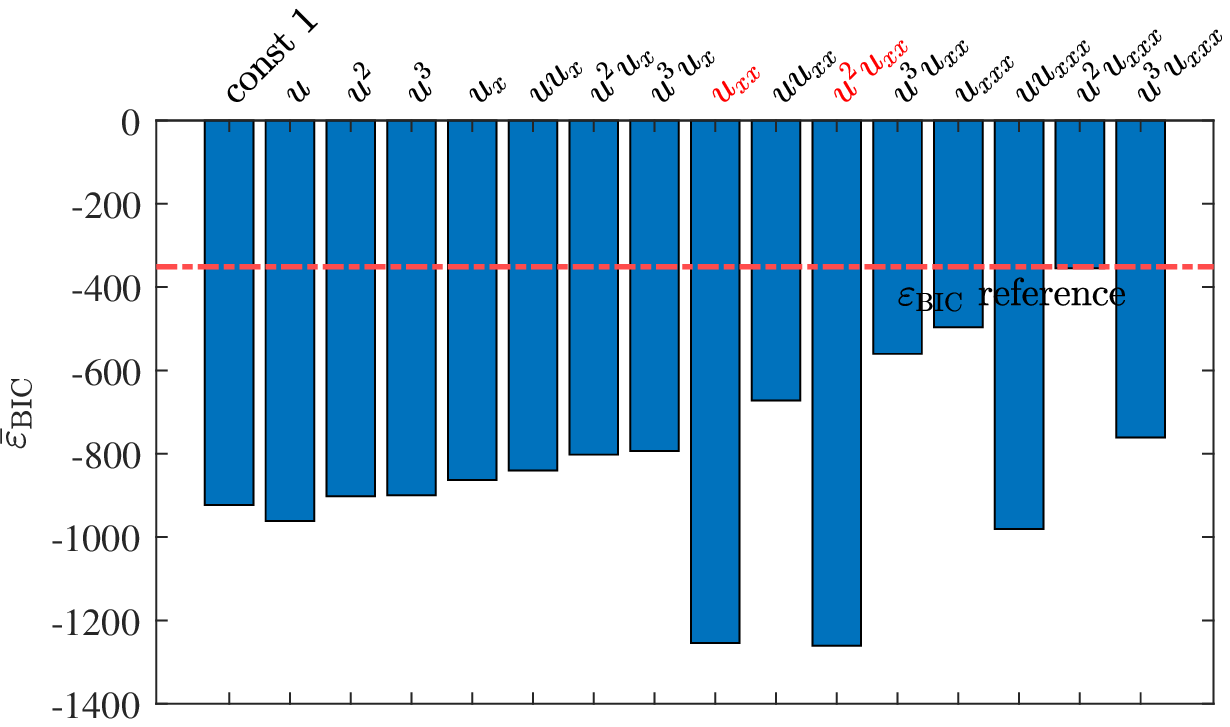}
	\includegraphics[scale=0.4]{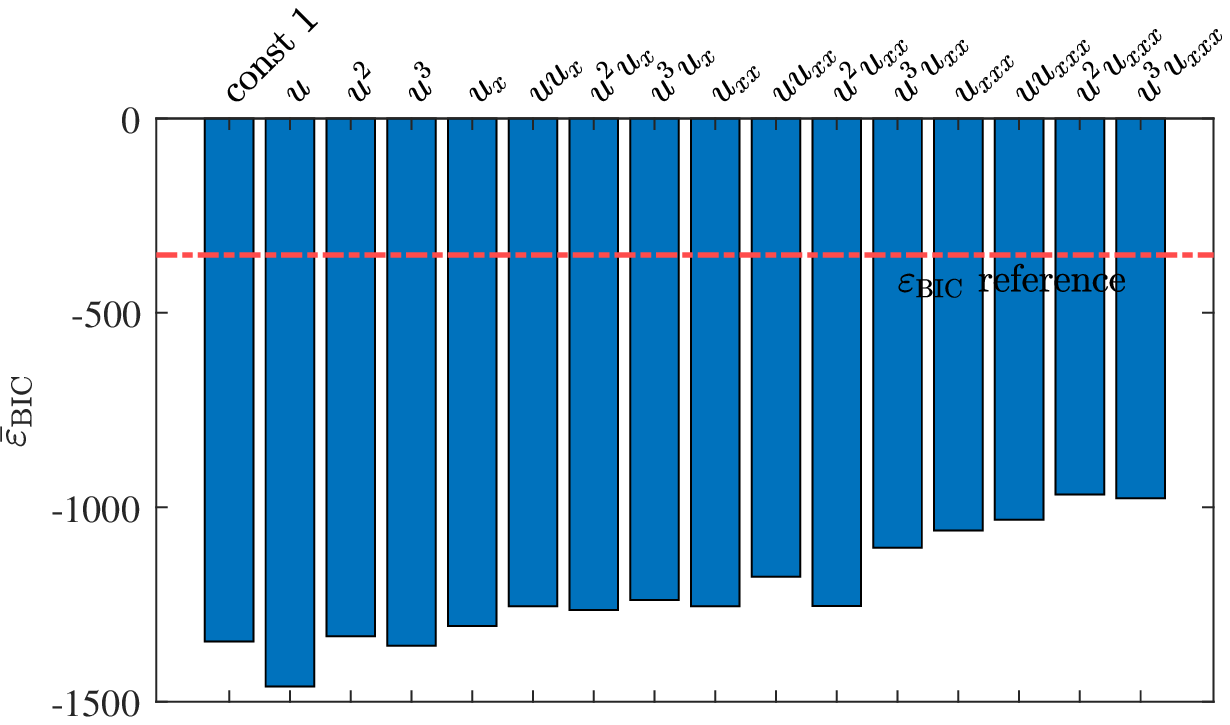}
	\caption{Mean BIC when dropping (a) or adding (b and c) a certain term in the library (Burgers equation: $u_t = -uu_x+\frac{0.01}{\pi}u_{xx}$).}
	\label{Figure:aveBIC} 
\end{figure}

Following the preprocessing steps that prepared the vector $\widetilde{\mathbf{U}}_{t}$ and the library matrix $\boldsymbol{\Theta}(\widetilde{\mathbf{U}})$, progressive sparse identification of the PDE is conducted using the $\psi$-algorithm presented in Algorithm \ref{alg1}. The form of the Burgers equation is identified progressively until when adding more terms to the learned PDE no longer considerably improves the regression accuracy. Without loss generality, the rest of this section uses the simulated clean data, and the results of cases with various levels of noise will be presented in section \ref{Sec:results} with other systems. First, the most contributive term(s) with index/indices $I_0$ are determined by comparing the increase of regression errors ($\varepsilon_{Reg}$ and $\varepsilon_{BIC}$) when dropping a certain term from the full list of terms in the library. Figures \ref{Figure:errReg} and \ref{Figure:errBIC} (a) plot the distributions of errors($\varepsilon_{Reg}$ and $\varepsilon_{BIC}$ respectively) when dropping each term in the library. It can be observed that deleting the $6^\mathrm{th}$ term $uu_x$ causes the largest regression error which is much larger than that when dropping any other term. This finding is further verified by comparing the mean errors as shown in Figures \ref{Figure:aveReg} and \ref{Figure:aveBIC} (a). Hence, the term $uu_x$ is first added to the PDE form with $I_0 = [6]$.

With $uu_x$ selected, the importance of other terms in the library are examined by adding each of them separately into the selected list and evaluating the improvement of regression. The results are shown in Figures \ref{Figure:errReg} and \ref{Figure:errBIC} (b) for the error distributions and Figures \ref{Figure:aveReg} and \ref{Figure:aveBIC} (b) for the mean errors. Comparison shows that adding the $9^\mathrm{th}$ term ($u_{xx}$) or the $11^\mathrm{th}$ term ($u^2u_{xx}$) can largely reduce the regression errors. Considering the complexity of the $11^\mathrm{th}$ term with higher nonlinearity at the same order of spatial derivative, the $9^\mathrm{th}$ term ($u_{xx}$) is added to the PDE and the selected index list is updated to $I_0 = [6,9]$. The $9^\mathrm{th}$ and $11^\mathrm{th}$ terms are not added simultaneously in this step since adding one of them may affect the importance of adding the other term regarding improving the regression accuracy. However, all possible outcomes of sparse regression with the $\psi$-algorithm will be investigated in the PDE solving/optimization step of the $\psi$-PDE method, which reevaluates the importance of each possible terms and finalizes the learned PDE of the observed system.

The same procedure is run to evaluate the importance of other terms, and the results (Figures \ref{Figure:errReg} to \ref{Figure:aveBIC} (c)) show that adding more terms to the PDE does not significantly improve the performance of regression but increases the complexity of the resulting model. Therefore, the $\psi$-algorithm for sparse regression will terminate with terms $uu_x$ and $u_{xx}$ selected to formulate the PDE for the given system, and the identified PDE in this step is $u_t = -0.9886uu_x+\frac{0.024}{\pi}u_{xx}$ with a least squares regression.

\begin{figure}[!h]
	\centering
	\tiny{(a) \hspace{130pt} (b)} \hspace{130pt} (c)\\
	\includegraphics[scale=0.35]{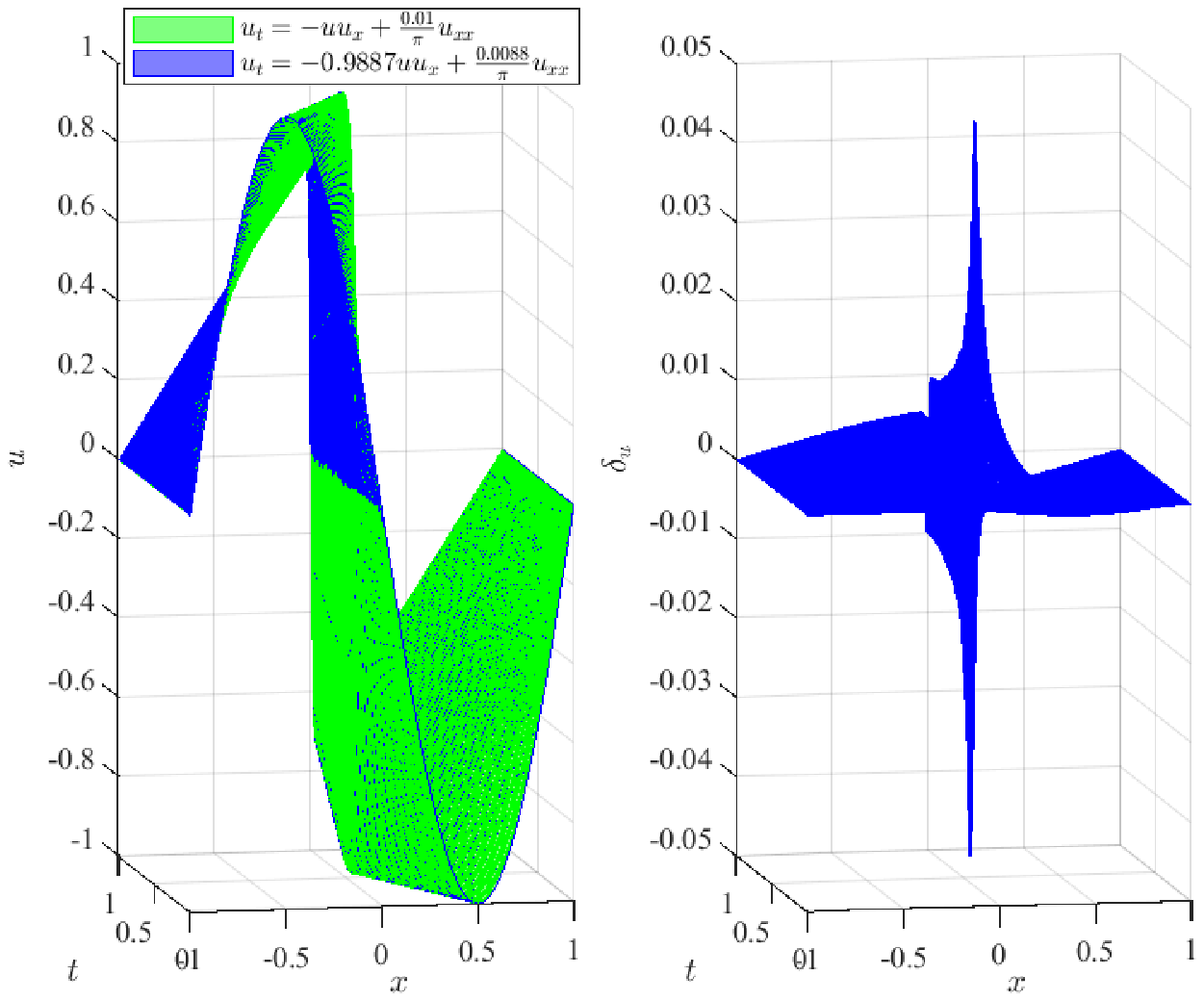}
	\includegraphics[scale=0.35]{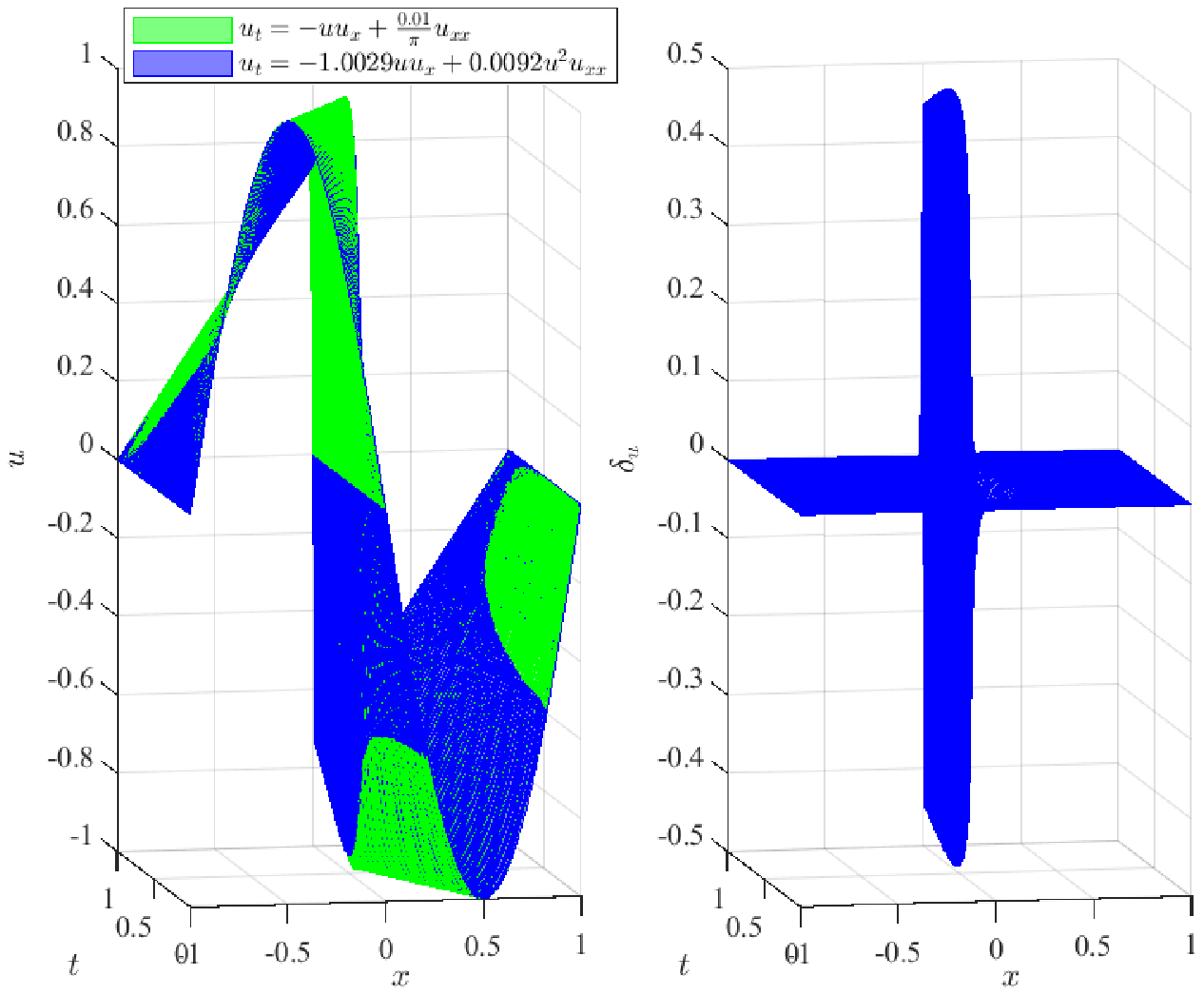}
	\includegraphics[scale=0.35]{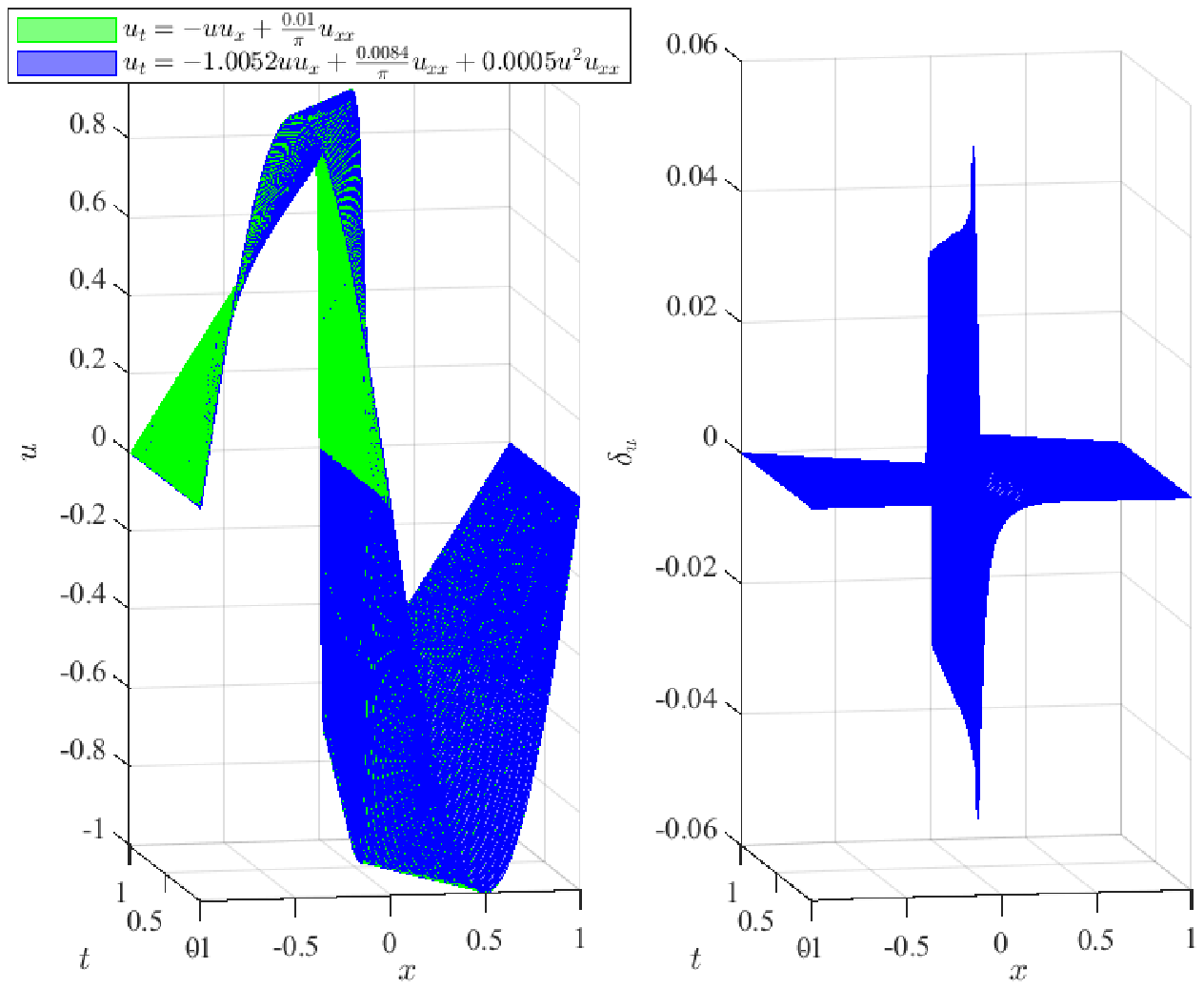}
	\caption{Comparison of all cadidate solutions with measured data (Burgers equation: $u_t = -uu_x+\frac{0.01}{\pi}u_{xx}$).}
	\label{Figure:solComp}
\end{figure}

With its form determined in sparse regression using the $\psi$-algorithm, the coefficients of each term of the identified PDE will be optimized in the last step of the $\psi$-PDE method by solving it with given/extracted initial/boundary conditions and comparing the solution with the measured data. The steepest descent method is used in this optimization. Figure \ref{Figure:solComp} (a) shows the resulting optimized PDE, its solution plotted together with the simulated data, and their difference $\delta_u$.The optimized PDE has exactly the same form with the ground truth PDE and accurate coefficients of both terms on the right side, and its solution is very close to the simulated data with the largest error below 5\%. Hence, the governing equation for the simulated system can be correctly identified using the proposed $\psi$-PDE method.

Considering that the $\psi$-algorithm does not necessarily yield a single solution especially in cases with noisy data, other possible solutions are analyzed to further investigate the importance of the selected/dropped terms as well as the robustness of the $\psi$-PDE method. In the sparse regression step, after selecting the $6^\mathrm{th}$ term $uu_x$, both the $9^\mathrm{th}$ term $u_{xx}$ and the $11^\mathrm{th}$ term $u^2u_{xx}$ are competitive candidates in terms of reducing the regression errors. Figures \ref{Figure:solComp} (b) and (c) shows the results of other possible solutions in the PDE solving/optimization step. In the second round in sparse regression, if the term $u^2u_{xx}$ is selected instead of $u_{xx}$, the $\psi$-algorithm will end with the PDE $u_t = -0.9886uu_x+0.0075u^2u_{xx}$. Figure \ref{Figure:solComp} (b) shows the results of optimization. It can be observed that the solution of the identified PDE has considerably large error at around $x=0$ (>50\%) though it matches well with the simulated data elsewhere. Additionally, one may suggest adding both  $u_{xx}$ and $u^2u_{xx}$ to the PDE, which yields the equation $u_t = -1.0101uu_x+\frac{0.0086}{\pi}u_{xx}+0.0047u^2u_{xx}$ in sparse regression. Figure \ref{Figure:solComp} (c) shows the results of optimization in this scenario. It shows that with the sacrifice of model parsimony, adding the term $u^2u_{xx}$ to the learned PDE does not largely reduce the difference of its solution from the simulated data. Moreover, one can observe that after optimization, the coefficient of the term $u^2u_{xx}$ becomes nearly ten times smaller, which further proves its insignificance in the governing equation of this system. This analysis with all candidate solutions of PDE learning ensures that the $\psi$-PDE method finally yields an equation that best captures the intrinsic underlying physics among all candidate solutions.

\section{Results and Discussions}\label{Sec:results}
This section presents and discusses the results of PDE learning using the $\psi$-PDE method with simulated noisy data from several canonical systems covering a number of scientific domains. 0 to 50\% noise is added to the numerically simulated clean data to demonstrate the effects of preprocessing and the robustness of the $\psi$-PDE method in cases with noisy measurements. Two 1D systems are investigated in section \ref{results1D}: (1) the dissipative system characterized by the 1D Burgers equation (as shown in section \ref{Sec:method})); (2) the traveling waves described by the Korteweg–de Vries (KdV) equation. To highlight the advantage of the $\psi$-PDE method in robustness, section \ref{results1D} first presents the results of learning the 1D Burgers equation using the PDE-FIND method developed in \cite{rudy2017data}. The lemma of hyperparameter tunning and the challenge in cases with noisy data in existing methods are demonstrated. Section \ref{results2D} presents the results of two 2D systems: (1) an extended dissipative system characterized by the 2D Burgers equation; (2) the lid-driven cavity flow corresponding to the 2D Navier Stokes equation. Due to the limited space in this article, the intermediate results are not presented (as done in section \ref{Sec:method}) for the cases studied in this section. Codes of demonstrated examples are available on the website: \href{https://github.com/ymlasu}{\textcolor{blue}{https://github.com/ymlasu}}.
\subsection{Discovering PDEs for 1D systems}\label{results1D}

Table \ref{Table:FIND} compares the results of learning the 1D Burgers equation ($u_t = -uu_x+\frac{0.01}{\pi}u_{xx}$) using the PDE-FIND method \cite{rudy2017data} considering various noise levels and different combinations of hyperparameters. The STRidge algorithm in the PDE-FIND method mainly has two hyperparameters, i.e., the regularizer $\lambda$ and the tolerance increment/ initial tolerance $d_{tol}$. From the authors' viewpoint, the hyperparameter tuning can be challenging in PDE learning where the ground truths and noise level are not known a priori. Therefore, a method may not be robust enough if the learning outcome is susceptible to the variation of hyperparameters. 

In cases with 0\% or 10\% noise, the PDE-FIND method yields identical PDE though with different hyperparameter combinations, as shown in the last column of Table \ref{Table:FIND}. Comparing the results with the ground truth PDE, one can find that correct terms are identified. However, without the PDE solving/optimization procedure, their coefficients cannot be optimized in this method. When the noise level increases to 20\%, the results become very sensitive to the variation of hyperparameters. With a different hyperparameter combination, the identified terms on the right hand side of PDE can be very different. Moreover, at this noise level, the terms of the ground truth PDE cannot be correctly identified though exhaustive trials of hyperparameter adjustment. Considering this lemma of hyperparameter tuning, especially in cases with a large level of noise, a robust method for PDE learning is needed when limited prior information is available for an unknown system.

\begin{table} [!h]
	\caption{Results of PDE learning using the PDE-FIND method in \protect\cite{rudy2017data} (Burgers equation: $u_t = -uu_x+\frac{0.01}{\pi}u_{xx}$).}
	\begin{tabular}{ |c|c|c|c| } 
		\hline
		Noise level & $\lambda$ & $d_{tol}$ & identified PDE \\ \hline
		0\% & -& - & $u_t = -0.7616uu_x+\frac{0.0114}{\pi}u_{xx}$\\ \hline
		10\% & -& - & $u_t = -0.6758uu_x+\frac{0.0099}{\pi}u_{xx}$\\ \hline
		20\% & $10^{-5}$ & 1.0 & $u_t = -0.0823uu_x-0.000001uu_{xxx}$\\ \hline
		20\% & $10^{-1}$ & 1.0 & $u_t = -0.772u^3u_{x} + 0.0031u^2u_{xx}$\\ \hline
		20\% & $10^{-5}$& 0.1 & $u_t = -0.5069uu_x+\frac{0.0034}{\pi}u_{xx}+0.0016u^2u_{xx} -0.00003uu_{xxx} + 0.00003u^3u_{xxx}$\\ \hline
		20\% & $10^{-1}$& 0.1 & $u_t = -0.772u^3u_{x} + 0.0031u^2u_{xx}$\\ \hline
	\end{tabular}\label{Table:FIND}
\end{table}

\begin{figure}[!h]
	\centering
	\tiny{(a) \hspace{130pt} (b)} \hspace{130pt} (c)\\
	\includegraphics[scale=0.35]{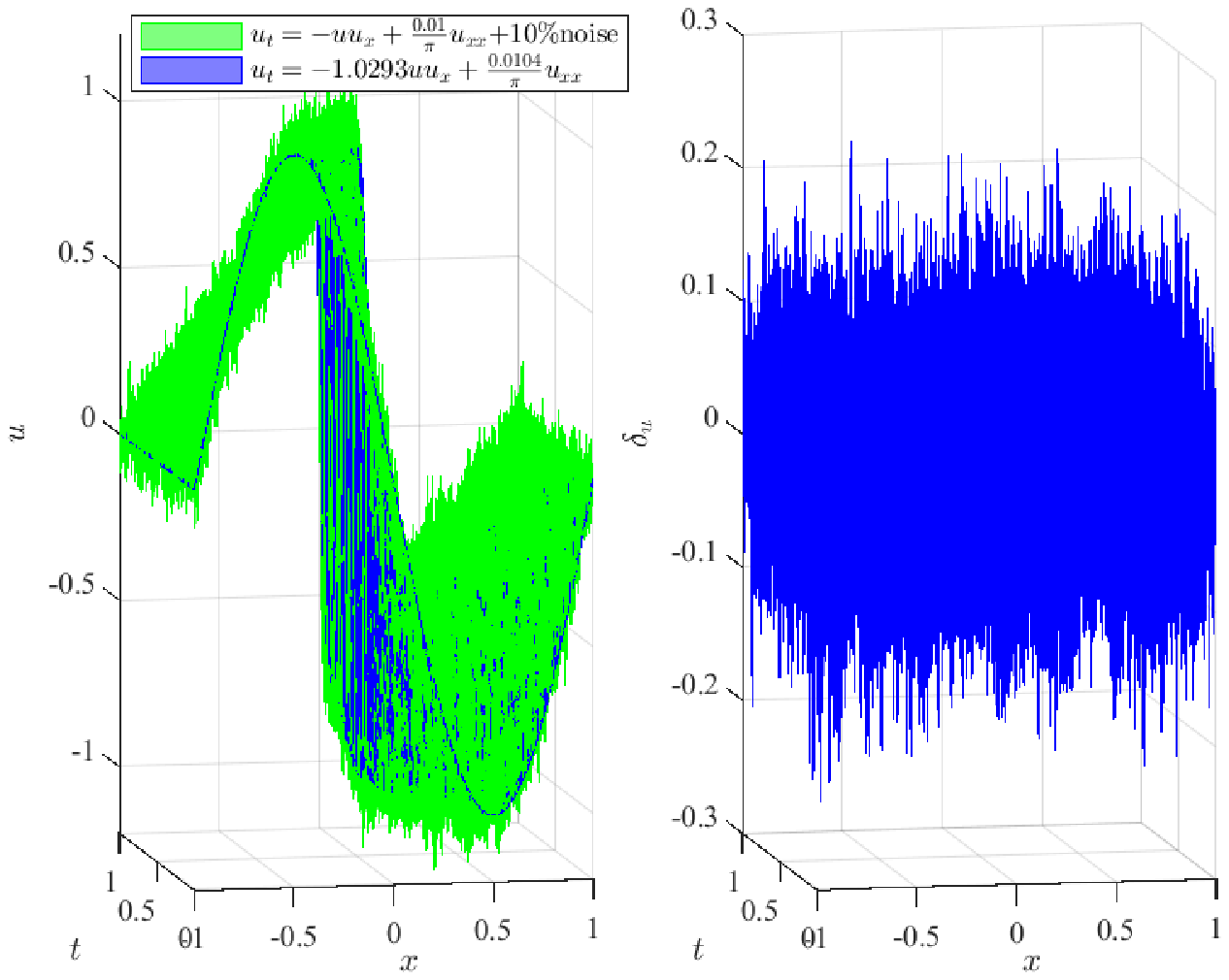}
	\includegraphics[scale=0.35]{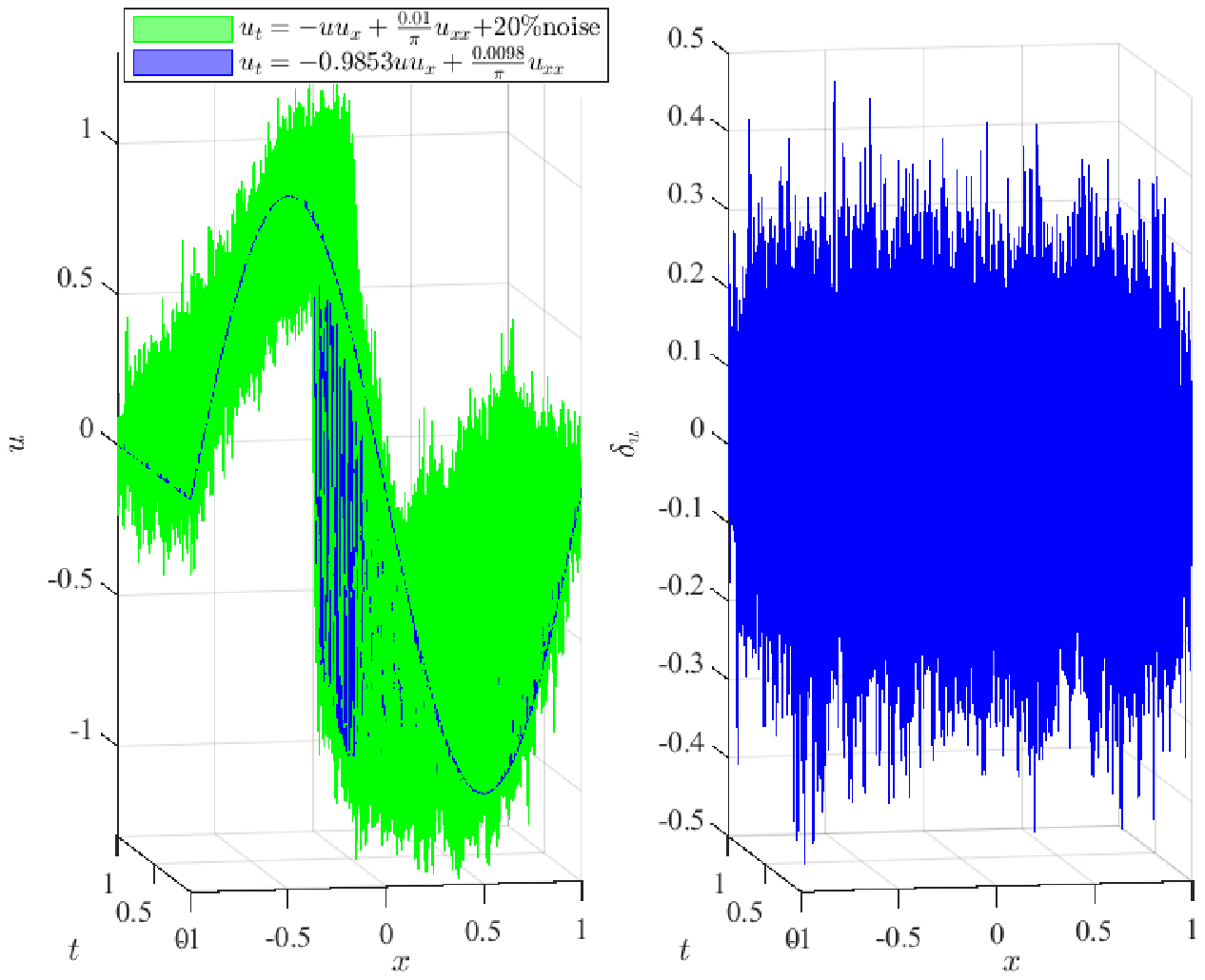}
	\includegraphics[scale=0.35]{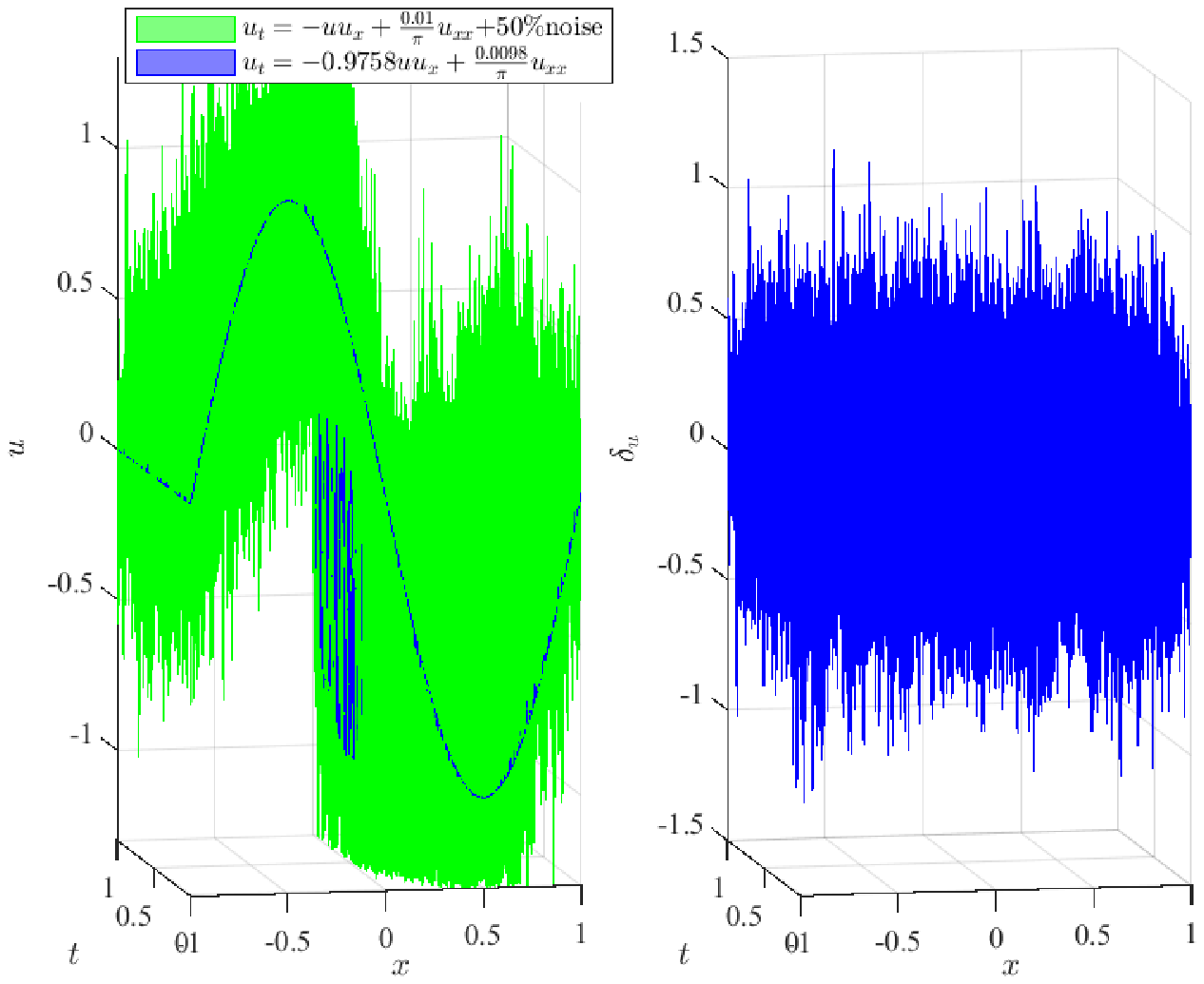}
	\caption{Results of PDE learning using the $\psi$-PDE method (Burgers equation: $u_t = -uu_x+\frac{0.01}{\pi}u_{xx}$). (a) 10\% noise; (b) 20\% noise; (c) 50\% noise.}
	\label{Figure:solNoisy}
\end{figure}

Table \ref{Table:burgersNoisy} lists the results of learning the Burgers equation ($u_t = -uu_x+\frac{0.01}{\pi}u_{xx}$) from noisy data using the proposed $\psi$-PDE method. Without complex and tricky hyperparameter setting/tuning, the $\psi$-PDE method yields identically correct PDE form even with data containing up to 50\% noise. Moreover, by virtue of the critical PDE solving/optimization step in the $\psi$-PDE method, the values of coefficients are very close to that of the ground truth PDE. It should be noted that the efficiency of the $\psi$-PDE method in learning correct PDEs holds beyond the noise levels investigated in this study. Figures \ref{Figure:solNoisy} (a) to (c) compare the solutions of learned PDEs with the measured noisy data. It can be observed that even in cases with significant noise, the $\psi$-PDE method never overfits the noise components in the measured data with redundant high-order nonlinear terms. Instead, this method always yields a governing equation that captures the most intrinsic invariants underlying the data.

\begin{table} [!h]
	\caption{Results of PDE learning using the $\psi$-PDE method (Burgers equation: $u_t = -uu_x+\frac{0.01}{\pi}u_{xx}$).}
	\begin{tabular}{ |c|c| } 
		\hline
		noise level &  identified PDE \\ \hline
		0\% & $u_t = -0.9887uu_x+\frac{0.0088}{\pi}u_{xx}$\\ \hline
		10\% &  $u_t = -1.0293uu_x+\frac{0.0104}{\pi}u_{xx}$\\ \hline
		20\% & $u_t = -0.9853uu_x+\frac{0.0098}{\pi}u_{xx}$\\ \hline
		50\% &  $u_t = -0.9758uu_x+\frac{0.0098}{\pi}u_{xx}$\\ \hline
	\end{tabular}\label{Table:burgersNoisy}
\end{table}

The second example of learning PDE from a 1D dynamical system examines the effectiveness of the $\psi$-method in correctly identifying governing equations containing higher-order spatial derivatives. This section considers a mathematical model of traveling waves on shallow water surfaces, i.e., the KdV equation with the form $u_t  =  \lambda_1uu_x + \lambda_2u_{xxx}$. The KdV equation can be used to characterize the evolution of many long 1D waves such as the ion acoustic waves in a plasma and acoustic waves on a crystal lattice \cite{raissi2019physics}. This study investigates the system described by the following KdV equation: $u_t = -uu_x-0.0025u_{xxx}$ with the initial condition $u(x,0)=\mathrm{cos}(\pi x)$ and periodic boundary conditions.
 
\begin{table} [!h]
	\caption{Results of PDE learning using the $\psi$-PDE method (KdV equation: $u_t = -uu_x-0.0025u_{xxx}$).}
	\begin{tabular}{ |c|c| } 
		\hline
		noise level &  identified PDE \\ \hline
		0\% & $u_t = -0.9998uu_x-0.0026u_{xxx}$\\ \hline
		10\% &  $u_t = -0.9757uu_x-0.0028u_{xxx}$\\ \hline
		20\% & $u_t = -0.9618uu_x-0.0031u_{xxx}$\\ \hline
		50\% &  $u_t = -1.0221.uu_x-0.0029u_{xxx}$\\ \hline
	\end{tabular}\label{Table:kdvNoisy}
\end{table}

\begin{figure}[!h]
	\centering
	\includegraphics[scale=0.80]{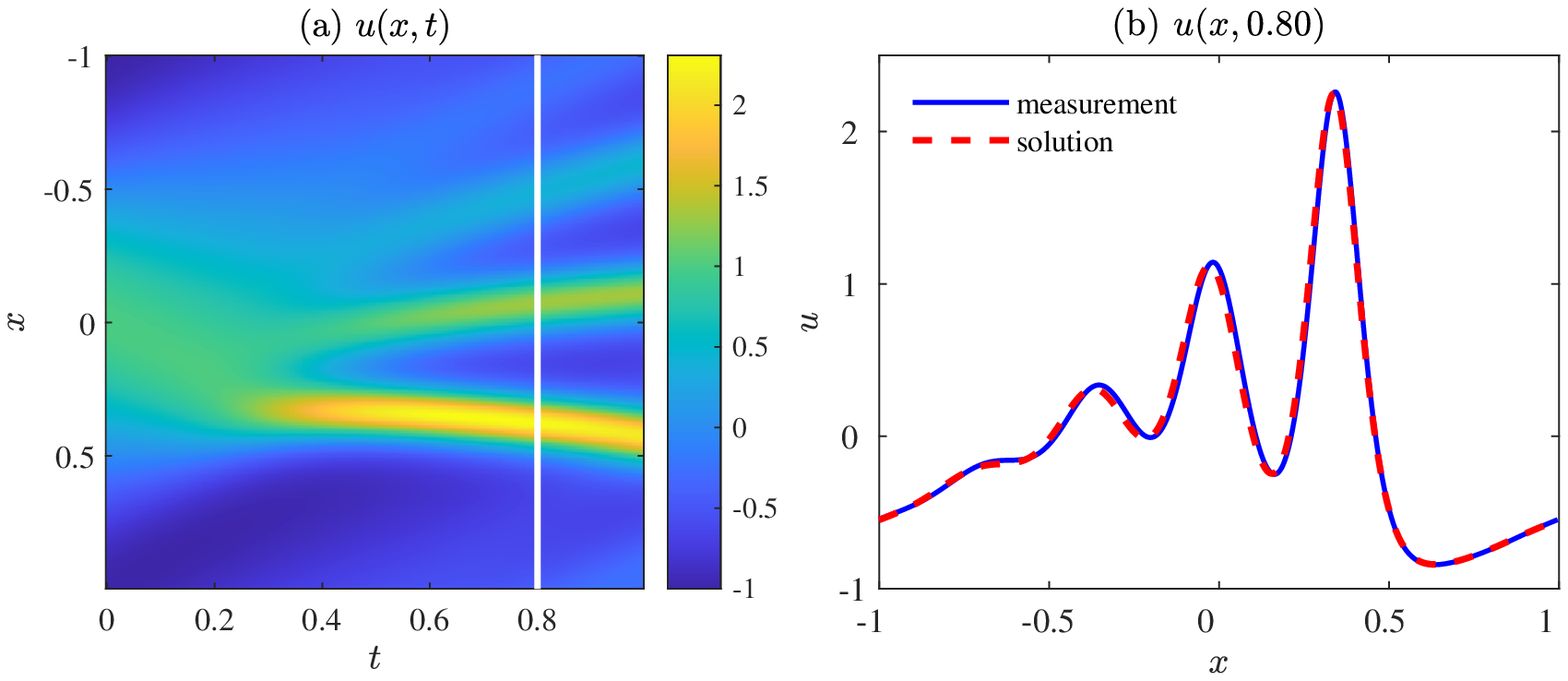}
	\caption{Results of learning the KdV equation ($u_t = -uu_x-0.0025u_{xxx}$: ) using the $\psi$-PDE method with clean data. The data/solution are plotted in 2D with values of $u$ denoted by color. (a) the measured data; (b) the measured $u$ with the solution to the learned PDE ($u_t = -0.9998uu_x-0.0026u_{xxx}$) at $t = 0.8$ sec.}
	\label{Figure:solKdVN0}
\end{figure}

\begin{figure}[!h]
	\centering
	\includegraphics[scale=0.80]{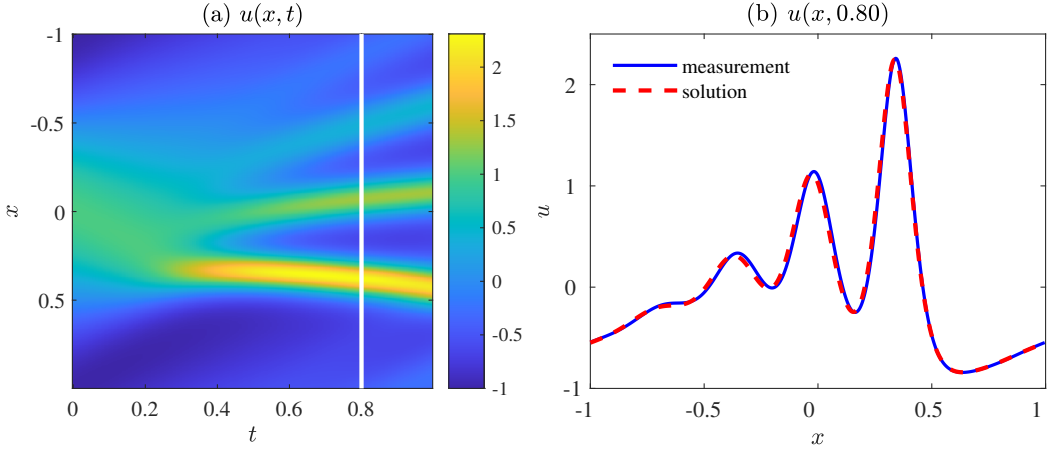}
	\caption{Results of learning the KdV equation ($u_t = -uu_x-0.0025u_{xxx}$) using the $\psi$-PDE method with data containing 50\% noise. (a) the measured data; (b) the measured $u$ with the ground truth (i.e., the solution to $u_t = -uu_x-0.0025u_{xxx}$) and the solution to the learned PDE ($u_t = -1.0221uu_x-0.0029u_{xxx}$) at $t = 0.8$ sec.}
	\label{Figure:solKdVN50}
\end{figure}

Table \ref{Table:kdvNoisy} summarizes the results of learning PDEs from the simulated traveling waves containing 0\% to 50\% noise. It shows that the correct PDE form with accurate coefficients can be identified using the $\psi$-PDE method in all cases. When implementing the $\psi$-PDE method, it was noted that the form of the KdV equation is easier to identify than that of the Burgers equation, because the $\psi$-algorithm for sparse regression does not yield more than one candidate PDE even for very noisy cases. The details can be found by running the code for this example available on the website: \href{https://github.com/ymlasu}{\textcolor{blue}{https://github.com/ymlasu}}. Figures \ref{Figure:solKdVN0} and \ref{Figure:solKdVN50} show the measured data and the its snapshot with the solution to the learned PDE at $t = 0.8$ sec. It can be found that even in case with 50\% noise, the $\psi$-PDE method can discover the intrinsic physics underlying the noisy data. Therefore, this example not only highlights the capability of the $\psi$-method of learning higher-order spatial derivatives but also further proves its the robustness in extracting the underlying physics from noisy measurements.

\subsection{Discovering PDEs for 2D systems}\label{results2D}
\begin{figure}[!h]
	\centering
	\includegraphics[scale=0.80]{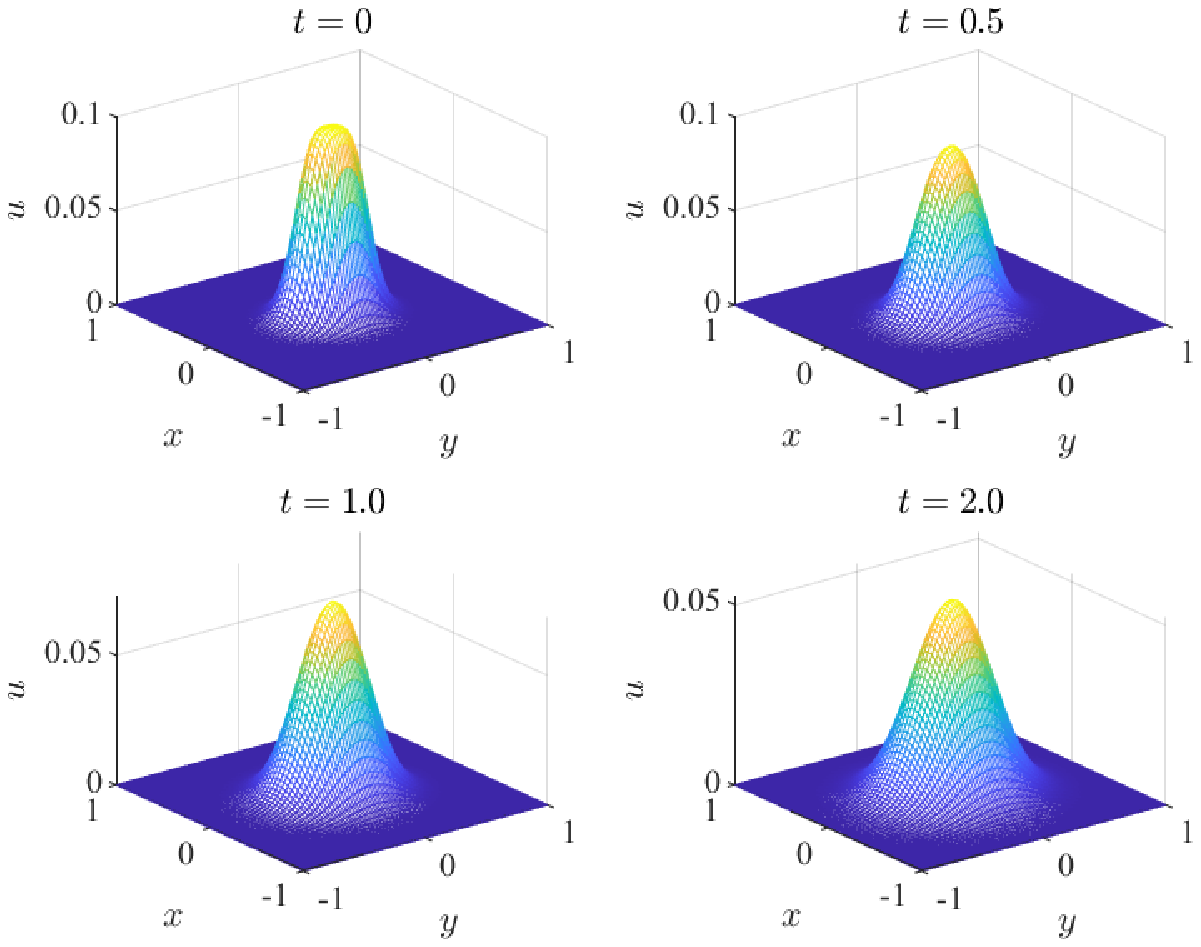}
	\caption{Dissipative system characterized by the 2D Burgers equation $u_t = -(uu_x+ uu_y)+0.01(u_{xx} +u_{yy})$.}
	\label{Figure:burgers2DN0}
\end{figure}

This section investigates the effectiveness of the $\psi$-PDE method in learning efficient governing equation(s) of 2D dynamical systems. First, the $\psi$-method is used to discover the physics of a 2D dissipative system characterized by the Burgers equation $u_t = -(uu_x+ uu_y)+0.01(u_{xx} +u_{yy})$ with the initial condition $u(x,y,0) = 0.1\mathrm{sech}(20x^2+25y^2)$ and periodic boundary conditions. Figure \ref{Figure:burgers2DN0} shows the snapshots of simulated data for this system. It should be noted that with the increase of dimensionality, the knowledge discovery of dynamical systems becomes more complex. Hence, in the sparse regression scheme of PDE learning, the library matrix $\boldsymbol{\Theta}$ should no longer be built in an exhaustive manner considering all possible combinations of polynomials to a certain power and spatial derivatives to a certain order, which will make the sparse regression problem intractable. Instead, $\boldsymbol{\Theta}$ is built with representative terms in multi-dimensional nonlinear dynamical systems (e.g., convective derivative $\boldsymbol{u} \cdot \nabla$, advective acceleration $(\boldsymbol{u} \cdot 	\nabla)\boldsymbol{u}$, and the Laplacian $	\nabla^2(\boldsymbol{u})$) and their products with polynomials. Table \ref{Table:2dBurgers} summarizes the results of learning PDEs from the 2D dissipative system using simulated clean and noisy data. It shows that the correct equation form with accurate coefficients can be successfully identified using the $\psi$-PDE method with data containing as much as 40\% random noise. When the noise level increases to 50\%, the $\psi$-algorithm fails to extract the most contributive terms from the library due to the severe influence of noise to the accuracy of calculated numerical derivatives. Advanced techniques of numerical differentiation such as automatic differentiation using graphical neural networks (GNN) \cite{baydin2017automatic} will be investigated in future studies.

\begin{table} [!h]
	\caption{Results of PDE learning using the $\psi$-PDE method (2D Burgers equation: $u_t = -(uu_x+ uu_y)+0.01(u_{xx} +u_{yy})$).}
	\begin{tabular}{ |c|c| } 
		\hline
		noise level &  identified PDE \\ \hline
		0\% & $u_t = -1.0028(uu_x+ uu_y)+0.0100(u_{xx} +u_{yy})$\\ \hline
		10\% &  $u_t = -1.0417(uu_x+ uu_y)+0.0100(u_{xx} +u_{yy})$\\ \hline
		20\% & $u_t = -1.0077(uu_x+ uu_y)+0.0100(u_{xx} +u_{yy})$\\ \hline
		30\% &  $u_t = -1.0723(uu_x+ uu_y)+0.0101(u_{xx} +u_{yy})$\\ \hline
		40\% &  $u_t = -1.0168(uu_x+ uu_y)+0.0101(u_{xx} +u_{yy})$\\ \hline
		50\% &  $-$\\ \hline
	\end{tabular}\label{Table:2dBurgers}
\end{table}

\begin{figure}[!h]
	\centering
	\includegraphics[scale=0.80]{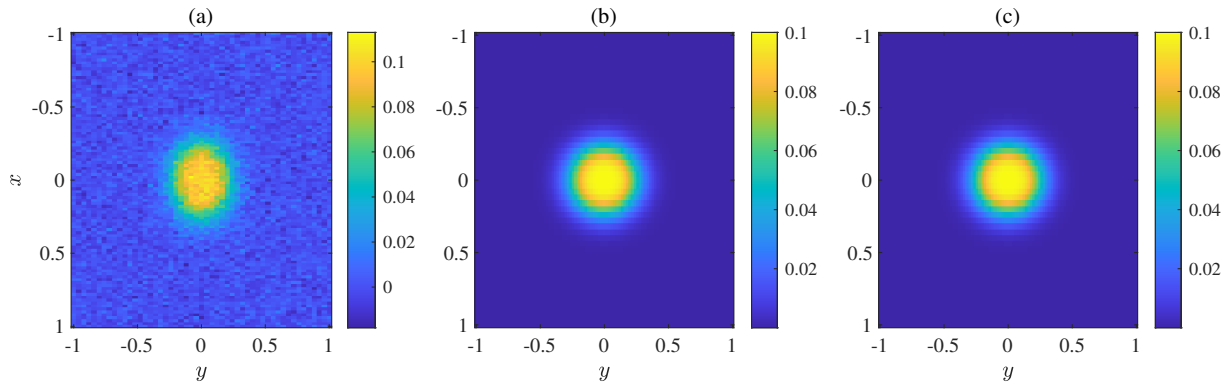}
	\caption{Results of learning the 2D Burgers equation ($u_t = -(uu_x+ uu_y)+0.01(u_{xx} +u_{yy})$) using the $\psi$-PDE method with data containing 40\% noise. The data/solution at $t = 0$ sec are plotted in 2D with values of $u$ denoted by color. (a) the measured data; (b) the ground truth (i.e., the solution to $u_t = -(uu_x+ uu_y)+0.01(u_{xx} +u_{yy})$); (c) the solution to the learned PDE ($u_t = -1.0168(uu_x+ uu_y)+0.0101(u_{xx} +u_{yy})$).}
	\label{Figure:solBurgers2dN40}
\end{figure}

Figure \ref{Figure:solBurgers2dN40} compares the solution to the learned PDE with the measured data containing 40\% noise and the solution to the ground truth PDE. One can observe that the proposed $\psi$-PDE method can discover the underlying truth behind the noisy measurements from the 2D dissipative system.

\begin{figure}[!h]
	\centering
	\includegraphics[trim=180 0 60 0,scale=0.45]{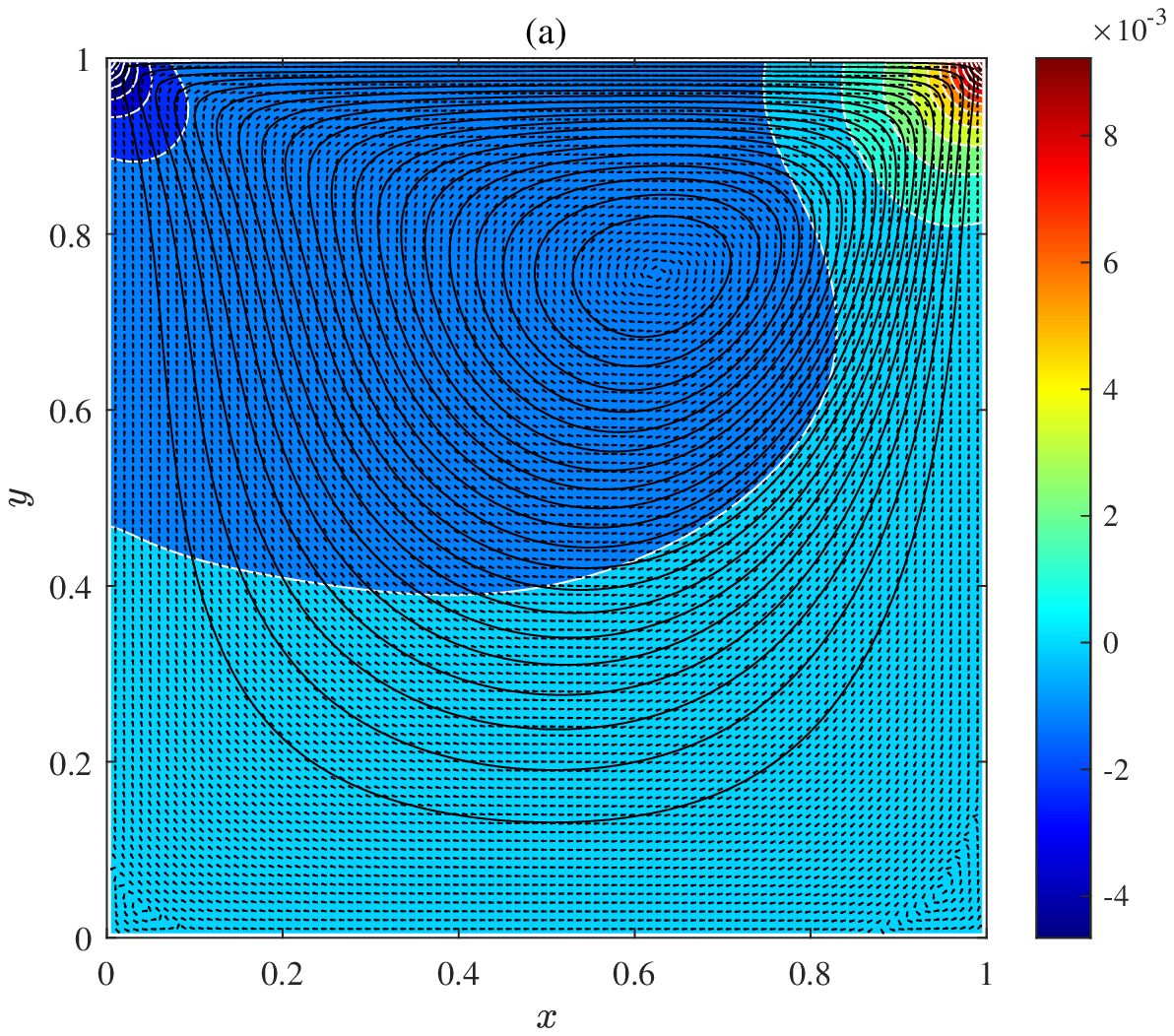}
	\includegraphics[trim=0 0 60 0,scale=0.45]{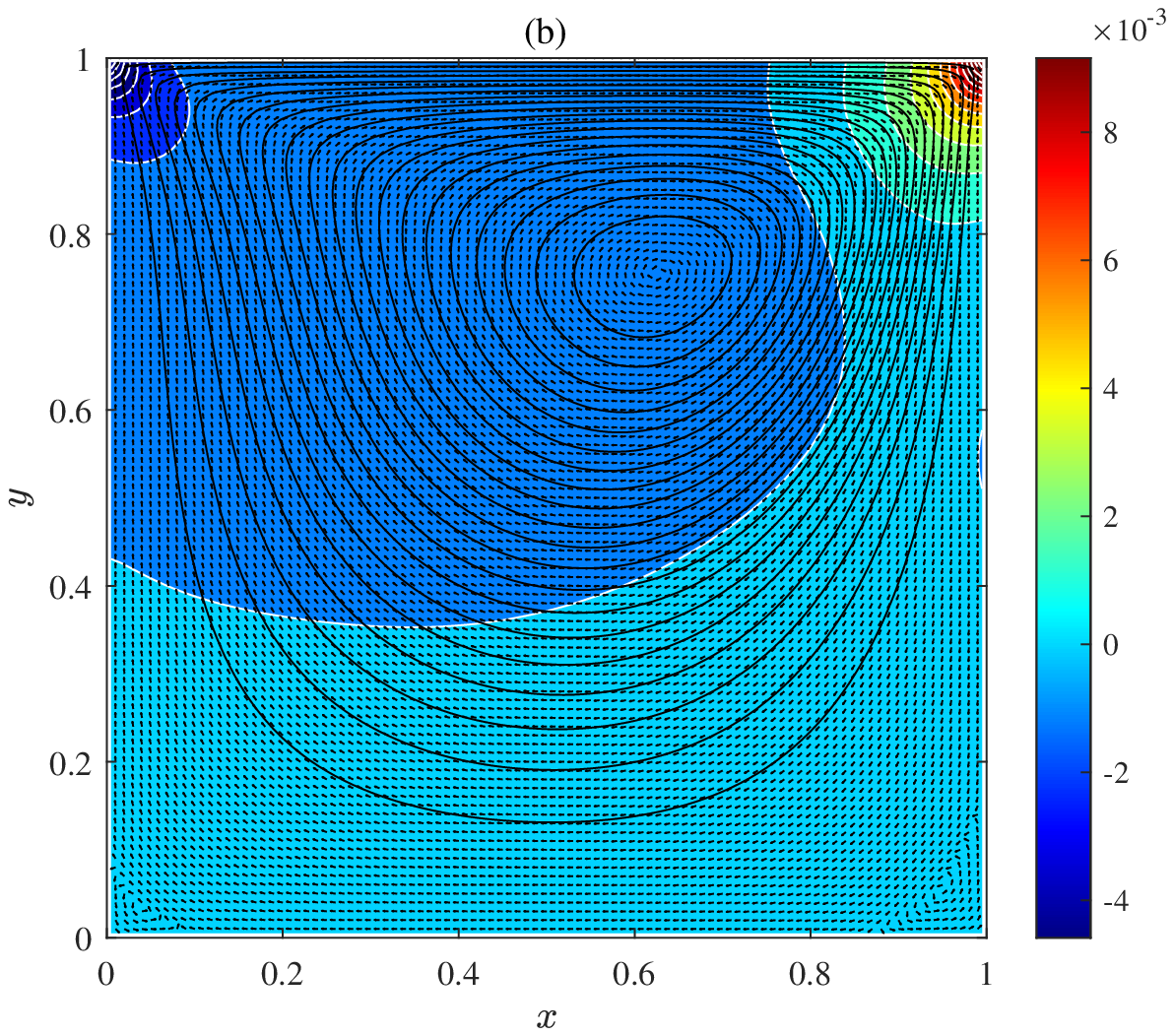}
	\includegraphics[trim=0 0 180 0,scale=0.45]{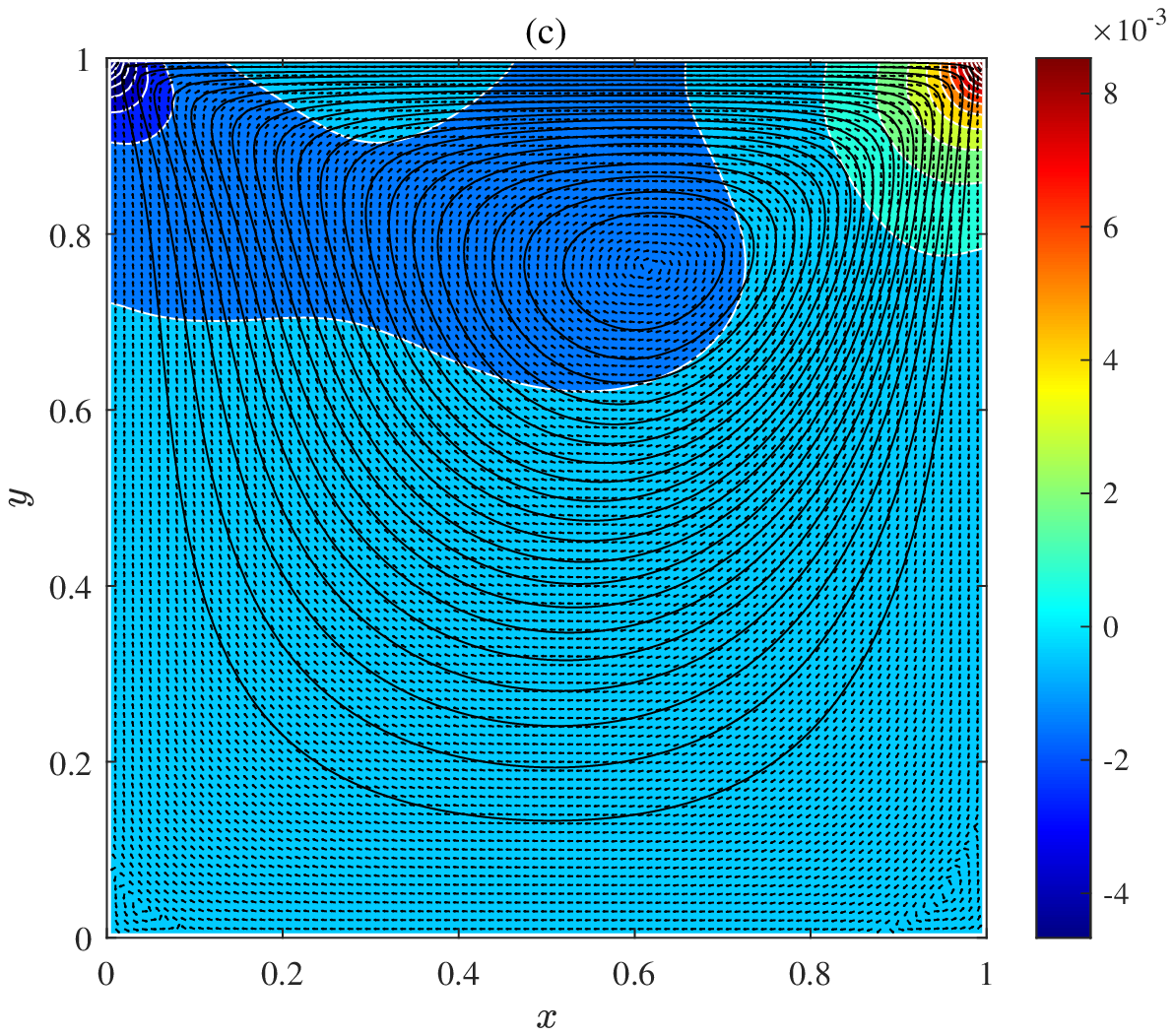}
	\caption{(a) Simulated lid driven flow (at $t = 4.0$ sec) characterized by the 2D Navier Stokes equation $\frac{\partial \boldsymbol{u}}{\partial t} = -(\boldsymbol{u}\cdot\nabla)\boldsymbol{u}-\nabla p + \frac{1}{Re}\nabla^2\boldsymbol{u}$ and $\nabla\cdot\boldsymbol{u}=0$ with $Re = 100$. (b) Solution to the learned PDE ($\frac{\partial \boldsymbol{u}}{\partial t} = -1.0060(\boldsymbol{u}\cdot\nabla)\boldsymbol{u}-1.0056\nabla p + \frac{1}{101.02}\nabla^2\boldsymbol{u}$) in case with 20\% noise. (c) Solution to the learned PDE ($\frac{\partial \boldsymbol{u}}{\partial t} = -0.8363(\boldsymbol{u}\cdot\nabla)\boldsymbol{u}-1.0069\nabla p + \frac{1}{101.91}\nabla^2\boldsymbol{u}$) in case with 50\% noise. The small arrows represent the velocity field; the color plot with contour lines denotes the pressure distribution; the closed contour lines are the streamlines.}
	
	\label{Figure:solNS}
\end{figure}

Another 2D system investigated in this study is the lid-driven cavity flow which is a benchmark problem for viscous incompressible fluid flow \cite{Zienkiewicz2005finite}. This study uses a geometry of a square cavity that is comprised of a lid on the top moving with a tangential unit velocity and three no-slip rigid walls. The velocity and pressure distributions are numerically simulated for a Reynolds number of 100. Figure \ref{Figure:solNS} (a) visualizes the solultion to this sytem at $t = 4.0$ sec, containing the velocity (small arrows) and pressure (color map with contour lines) distributions as well as the streamlines.

Table \ref{Table:NS} lists the results of PDE learning using the $\psi$-PDE method with simulated data containing various levels of noise. The correct PDE form can be learned for cases with as much as 50\% noise. Figure \ref{Figure:solNS} (b) illustrates the solution to the learned PDE in case with 20\% noise, which is almost identical to the solution to the ground truth PDE (Figure \ref{Figure:solNS} (a)). This observation further validates the effectiveness of the $\psi$-PDE method in distilling the underlying intrinsic physics. When the noise level increases to 50\%, the coefficient of the advective acceleration term (i.e.,($\boldsymbol{u}\cdot\nabla)\boldsymbol{u}$) is about 14\% different from the real value, as shown in the bottom row of Table \ref{Table:NS}. As a result, its solution especially the pressure distribution is different from that of the ground truth PDE, as shown in Figure \ref{Figure:solNS} (c). This difference demonstrates the challenge of PDE learning for a multidimensional system especially in the presence of large noise, which will be further investigated in the future work.

\begin{table} [!h]
	\caption{Results of PDE learning using the $\psi$-PDE method (2D Navier Stokes equation: $\frac{\partial \boldsymbol{u}}{\partial t} = -(\boldsymbol{u}\cdot\nabla)\boldsymbol{u}-\nabla p + \frac{1}{100}\nabla^2\boldsymbol{u}$ with $\nabla\cdot\boldsymbol{u}=0$).}
	\begin{tabular}{ |c|c| } 
		\hline
		noise level &  identified PDE \\ \hline
		0\% & $\frac{\partial \boldsymbol{u}}{\partial t} = -0.9649(\boldsymbol{u}\cdot\nabla)\boldsymbol{u}-0.9920\nabla p + \frac{1}{104.50}\nabla^2\boldsymbol{u}$\\ \hline
		10\% & $\frac{\partial \boldsymbol{u}}{\partial t} = -0.9678(\boldsymbol{u}\cdot\nabla)\boldsymbol{u}-0.9769\nabla p + \frac{1}{104.08}\nabla^2\boldsymbol{u}$ \\ \hline
		20\% & $\frac{\partial \boldsymbol{u}}{\partial t} = -1.0060(\boldsymbol{u}\cdot\nabla)\boldsymbol{u}-1.0056\nabla p + \frac{1}{101.02}\nabla^2\boldsymbol{u}$\\ \hline
		50\% & $\frac{\partial \boldsymbol{u}}{\partial t} = -0.8363(\boldsymbol{u}\cdot\nabla)\boldsymbol{u}-1.0069\nabla p + \frac{1}{101.91}\nabla^2\boldsymbol{u}$\\ \hline
	\end{tabular}\label{Table:NS}
\end{table}

\section{Summary and Further Discussions} \label{Section:Conclusion}
In this study, a robust data-driven method (i.e., the $\psi$-PDE method) is proposed for discovering the underlying physics of a given system from measured data. Investigating and improving its robustness is critical for effectively distilling the intrinsic law underlying a complex novel dynamical system. The $\psi$-PDE method has been tested on various systems in sufficiently challenging scenarios regarding the model complexity and noise intensity, and the identification results approve its effectiveness and generality. Compared with the state-of-the-art methods in the literature, the $\psi$-PDE method is advantageous in its robustness since it requires the least effort in hyperparameter tuning which should be obviated in identifying the governing equation(s) of an unknown system.

This work was inspired by the pioneering work presented in \cite{schmidt2009distilling} which discourages automatic discovery of natural laws. After examining the challenge of automatic identification via sparse regression as presented in \cite{rudy2017data} and its following works (such as \cite{berg2019data} and \cite{chen2020deep}), this study borrows the ``nonautomatic'' idea in \cite{schmidt2009distilling} and yields more than one candidate solutions through sparse regression with the $\psi$-algorithm; finally, the representativeness of all candidate models are evaluated through solving and optimizing their respective PDE forms taking the measured data as the ground truth and objective. It has been demonstrated that the $\psi$-PDE method always yields equations that capture the most intrinsic physics of the observed system. In comparison, most existing methods yield a unique PDE for a given system and do not allow evaluating its effectiveness in characterizing the system or optimizing its representativeness.

While it increases the robustness of PDE learning, the PDE solving/optimization step in the $\psi$-PDE method considerably increases the computational cost especially for a high-dimensional system. However, this challenge can be solved by parallel computing and/or surrogate modeling, which will be investigated in future studies for more complex systems. The future work may also include applying the $\psi$-PDE method real operating dynamical systems to further examine or improve its effectiveness in knowledge discovery.

\section{Acknowledgments}
The research reported in this paper was supported by funds from NASA University Leadership Initiative Program (Contract No. NNX17AJ86A, Project Officer: Dr. Anupa Bajwa, Principal Investigator: Dr. Yongming Liu). The support is gratefully acknowledged.

\section*{References}

\bibliographystyle{elsarticle-num}
\bibliography{PDE_Learn}

\end{document}